\documentclass[10pt]{article}
\usepackage{amssymb,amsbsy,amsmath,amsfonts,amssymb,amscd}
\usepackage[english]{babel}
\usepackage[T1]{fontenc}
\usepackage{indentfirst}
  \def\t{{\mathrm{t}}}

  \def\CC{{\mathbb{C}}}
\def\DD{{\mathbb{D}}}

 \def\NN{{\mathbb{N}}} 
 \def\QQ{{\mathbb{Q}}} \def\RR{{\mathbb{R}}}
 \def\TT{{\mathbb{T}}} 
  
 \def\ZZ{{\mathbb{Z}}}
\def\cA{{\mathcal{A}}}  \def\cC{{\mathcal{C}}}
\def\cD{{\mathcal{D}}}  
 \def\cH{{\mathcal{H}}}

  \def\cR{{\mathcal{R}}}

 \let\leq=\leqslant
 \let\geq=\geqslant
\def\og{\leavevmode\raise.3ex\hbox{$\scriptscriptstyle\langle\!\langle$~}}
\def\fg{\leavevmode\raise.3ex\hbox{~$\!\scriptscriptstyle\,\rangle\!\rangle$}}
\newtheorem{statement}{}
\newtheorem{theoreme}[statement]{Theorem}
\newtheorem{lemme}[statement]{Lemma}

\newtheorem{proposition}[statement]{Proposition}

\newtheorem{corollaire}[statement]{Corollary}

\title{Composition operators on the Wiener-Dirichlet algebra}
\author{Fr\'ed\'eric Bayart, Catherine Finet, Daniel Li, Herv\'e Queff\'elec}

\date{\footnotesize \today}

\begin{document}

\maketitle

\emph{{\bf Abstract.} We study the composition operators on an algebra of 
Dirichlet series, the analogue of the Wiener algebra of absolutely convergent Taylor series, 
which we call the Wiener-Dirichlet algebra. The central issue is to understand the connection between the 
properties of the operator and of its symbol, with special emphasis on the compact, automorphic, or 
isometric character of this operator. We are led to the intermediate study of algebras of 
functions of several, or countably many, complex variables.}
\bigskip

\noindent{2000 Mathematics Subject Classification -- 
primary: 47 B 33, secondary: 30 B 50 -- 42 B 35
}

\section{Introduction}

Let $A^+ = A^+(\TT)$ be the Wiener algebra of absolutely convergent
Taylor series in one variable~: $f \in A^+$ if and only if
\begin{displaymath}
f(z) = \sum_{n=0}^\infty a_n z^n,\hskip 3mm \mbox{with}\hskip 2mm 
\|f\|_{A^+} = \sum_{n=0}^\infty |a_n| < +\infty.
\end{displaymath}
It is well-known that $A^+$ is a commutative, unital Banach
algebra with spectrum $\overline \DD$, the closed unit disk.  If
$\phi : \DD \to \overline \DD$ is analytic, the composition operator
$C_\phi$ with symbol $\phi$ is formally defined by $C_\phi(f) = f\circ \phi$. \par 
\indent
Newman \cite{Ne} studied those symbols $\phi$ generating
bounded composition operators $C_\phi \colon A^+ \to A^+$, and proved in particular
the following:
\begin{enumerate}
\renewcommand{\labelenumi}{(\alph{enumi})}
\item $C_\phi$ maps $A^+$ into itself if and only if $\phi \in A^+$
and $\|\phi^n\|_{A^+} = O\,(1)$ as $n \to \infty$ (\textit{e.g.} 
$\phi(z) = 5^{-1/2}(1+z-z^2)$): this happens if and only if all maximum points 
$\theta_0$ of $|\phi ({\rm e}^{i\theta})|$ are ``\emph{ordinary points}'', \emph{i.e.} 
if and only if we have, as $t\to 0$: 
\begin{displaymath}
\log \phi({\rm e}^{i(\theta_0 +t)})=\alpha_0+\alpha_1 t+\alpha_k t^k+\cdots,
\end{displaymath}
where $k>1$ and $\alpha_k\neq 0$ is not purely imaginary;
\item if moreover $|\phi(e^{it})| = 1$, one must have
$\phi(z) = az^d$, with $|a|=1$ and $d\in \NN$;
\item $C_\phi \colon A^+ \to A^+$ is an automorphism if and only if
$\phi(z) = az$, with $|a|=1$.
\end{enumerate}
Harzallah (see \cite{K}) also proved
that:
\begin{enumerate}
\renewcommand{\labelenumi}{(\alph{enumi})}
\setcounter{enumi}{3}
\item $C_\phi \colon A^+ \to A^+$ is an isometry if and only if
$\phi(z) = az^d$, with $|a|=1$ and $d \in \NN$.
\end{enumerate}
\par
The aim of this paper is to perform a similar study for the
``\emph{Wiener-Dirichlet}'' algebra $\cA^+$ of absolutely convergent
Dirichlet series: $f \in \cA^+$ if and only if 
\begin{displaymath}
f(s) =\sum_{n=1}^\infty a_n n^{-s},\hskip 3mm \hbox{with}\hskip 2mm 
\|f\|_{\cA^+} = \sum_{n=1}^\infty |a_n| < +\infty.
\end{displaymath} 
\par
\noindent $\cA^+$ is a commutative, unital Banach algebra, with the
following multiplication (quite different from the one for
Taylor series):
\begin{displaymath}
\left( \sum_{n=1}^\infty a_n n^{-s} \right) 
\left( \sum_{n=1}^\infty b_n n^{-s} \right) = 
\sum_{n=1}^\infty c_n n^{-s}, \mbox{~with~}
c_n = \sum_{ij=n} a_i b_j.
\end{displaymath}
$\cA^+$ can also be interpreted as a space of analytic functions on $\CC_0$
(where in general we denote by $\CC_\theta$ the vertical half-plane 
${\cR e}\,s > \theta$). The study of function spaces formed by
Dirichlet series has gained some recent interest (see the papers of
Hedenmalm-Lindqvist-Seip \cite{HLS}, Gordon-Hedenmalm \cite{GH},
Bayart \cite{B1}, \cite{B2}, Finet-Queff\'elec-Volberg \cite{FQV},
Finet-Queff\'elec \cite{FQ}, Finet-Li-Queff\'elec \cite{FLQ}, 
Mc Carthy \cite{McC}).  Now, a method due to Bohr (see for example \cite{Q}) 
identifies the algebra $\cA^+$ with the algebra $A^+(\TT^\infty)$ formed by the 
absolutely convergent Taylor series in countably many variables (this point of view, 
which allows to identify the spectrum of $\cA^+$ as $\overline \DD^\infty$, the 
spectrum of $A^+(\TT^\infty)$, has been used by Hewitt and
Williamson \cite{HeW}, among others, to prove the following Wiener type tauberian
Theorem~: ``{\sl If $f \in \cA^+$ and $|f(s)| \geq \delta > 0$ for
$s \in \CC_0$, then $1/f \in \cA^+$}'').\par
Let us recall the way this identification is carried out. 
Let $(p_j)_{j \geq 1}$ be the increasing sequence of prime numbers 
($p_1 = 2$, $p_2 = 3$, $p_3=5$, \ldots). If
\begin{displaymath}
f(z) = \sum_{\alpha \in \NN_0^{(\infty)}} a_\alpha z^\alpha
\mbox{~~with~~} \|f\|_{A^+(\TT^\infty)} = \sum_\alpha |a_\alpha| < +\infty,
\end{displaymath}
where, as usual, we set $\alpha = (\alpha_1, \ldots,\alpha_r, 0,0, \ldots)$ and  
$z^\alpha = z_1^{\alpha_1} \ldots z_r^{\alpha_r}$ for $z = (z_j)_{j \geq 1}$,  
then $\Delta \colon \cA^+ \to A^+(\TT^\infty)$ is defined by: 
\begin{displaymath}
\Delta \left( \sum_{n=1}^\infty a_n n^{-s} \right) = 
\sum_{n=1}^\infty a_n z_1^{\alpha_1} \ldots z_r^{\alpha_r},
\end{displaymath}
if $n = p_1^{\alpha_1}\ldots p_r^{\alpha_r}$ is the decomposition of $n$ in prime
factors. $\Delta$ is an isometric isomorphism. Moreover, we shall need two 
more facts about $\Delta$. For
$s \in \CC_0$, we set $z^{[s]} = (p_j^{-s})_j$. We then have:
\begin{eqnarray}
\Delta f (z^{[s]}) &=& f(s), \mbox{~for any $f \in \cA^+$ and
any $s \in \CC_0$} \label{equation32}\\ 
\|\Delta f\|_\infty &=& \|f\|_\infty \mbox{~for each $f \in \cA^+$},\hfill
\label{equation33}
\end{eqnarray}
where we set $\|f\|_\infty = \sup\limits_{s \in \CC_0} |f(s)|$ and
$\|\Delta f\|_\infty = \sup\limits_{z \in {\textbf B}} |\Delta f(z)|$, with 
${\textbf B}=\{z=(z_j)_{j\geq 1}\in \DD^\infty\,;\ 
z_j\mathop{\longrightarrow}\limits_{j\to+\infty} 0\}$.
Indeed, if $f(s) = \sum_1^\infty a_n n^{-s}$, we have:
\begin{displaymath}
\Delta f(z^{[s]}) = \sum_{n=1}^\infty a_n(p_1^{-s})^{\alpha_1} \ldots
(p_r^{-s})^{\alpha_r} =
\sum_{n=1}^\infty a_n (p_1^{\alpha_1} \ldots p_r^{\alpha_r} )^{-s} =
f(s).
\end{displaymath}
On the other hand, let $z = (z_j)_{j \geq 1} \in {\textbf B}$.  Fix an
integer $N$, let $k = \pi(N)$ be the number of primes not
exceeding $N$, and 
$S_N(z) = \sum_{n=1}^N a_n z_1^{\alpha_1}\ldots z_k^{\alpha_k}$, 
with $n = p_1^{\alpha_1} \ldots p^{\alpha_k}_k$. Pick $\sigma > 0$ such 
that $|z_j| \leq p_j^{-\sigma}$, $1 \leq j \leq k$.  Due to the rational
independence of $\log p_1, \ldots, \log p_k$ and to the Kronecker  
Approximation Theorem (\cite{HST}, Corollary 4, page 23), the points $(p_j^{-it})_{1 \leq j \leq k}$, 
$t \in \RR$, are dense in the torus $\TT^k$, so that the
maximum modulus principle for the polydisk $\DD^k$ gives:
\begin{eqnarray*}
|S_N(z)| & \leq & \sup_{|w_j|=p_j^{-\sigma}} 
\Big|\sum_{n=1}^N a_n w_1^{\alpha_1} \ldots w_k^{\alpha_k} \Big| =
\sup_{{\cR e}\,s = \sigma} \Big| 
\sum_{n=1}^N a_n (p_1^{-s})^{\alpha_1}\ldots (p_k^{-s})^{\alpha_k} \Big| \\
& = & \sup_{{\cR e}\,s = \sigma} \Big| 
\sum_{n=1}^N a_n n^{-s} \Big|
\leq \Big\| \sum_{n=1}^N a_n n^{-s} \Big\|_\infty.
\end{eqnarray*}
Hence $\|S_N\|_\infty \leq \| \sum_1^N a_n n^{-s}\|_\infty$. Letting $N$ 
tend to infinity gives $\|\Delta f\|_\infty \leq \|f\|_\infty$, which 
proves (\ref{equation33}), 
since we trivially have $\|\Delta f\|_\infty \geq \|f\|_\infty$.\par\smallskip

\indent In this paper, we use the identification proposed above to obtain 
results similar to (a), (b), (c) and (d) for $\cA^+$. This leads to an 
intermediate study of composition operators on the
algebras $A^+(\TT^\infty)$ and $A^+(\TT^k)$ (the $k$-dimensional
analog of $A^+(\TT^\infty)$). Accordingly, the paper is organized as
follows:\par
\noindent In Section 2, we give necessary as well as sufficient 
conditions for boundedness and compactness of 
$C_\phi : \cA^+ \to \cA^+$, and study in detail some specific examples. In Section 3, 
we study the automorphisms of the algebras $A^+(\TT^k)$, $A^+(\TT^\infty)$,
$\cA^+$. In Section 4, we study the isometries of those algebras, and we
point out some specific differences between the finite and
infinite-dimensional cases. Section 5 is devoted to some concluding remarks 
and questions. \par\smallskip

A word on the definitions and notations: we will say that
integers $2 \leq q_1 < q_2 < \ldots$ are multiplicatively
independent if their logarithms are rationally independent in the
real numbers; equivalently, if any integer $n \geq 2$ can be
expressed as $n = q_1^{\alpha_1} \ldots q_r^{\alpha_r}$, $\alpha_j
\in \NN_0$, in at most one way (e.g. $q_1 = 2$, $q_2 = 6$, $q_3 =
30$).  We shall denote by ${\cD}$ the space of functions 
$\varphi \colon \CC_0 \to \CC$ which are analytic, and moreover representable as a
convergent Dirichlet series $\sum_1^\infty c_n n^{-s}$ for
${\cR e}\,s$ large enough. $\cD$ is also called the space of convergent
Dirichlet series. For example, if 
$\psi(s) = (1-2^{1-s})\zeta(s)$, and $\varphi(s) = \psi(s-a)$, $\varphi$ is entire, and
representable as $\sum_1^\infty (-1)^{n-1} n^a n^{-s}$ for ${\cR e}\,s > a$.  
$\TT$ denotes the unit circle, and plays no role in
the definition of $A^+(\TT^k)$ and $A^+(\TT^\infty)$, although $\TT^k$
(resp. $\TT^\infty$) might be viewed as the {\v S}hilov boundary of
$A^+(\TT^k)$ (resp. $A^+(\TT^\infty)$). As usual, we set $\NN = \{1,2,\ldots\}$ and
$\NN_0 = \{0,1,2,\ldots\} = \NN \cup \{0\}$. Recall that $\CC_\theta$ is the 
vertical half-plane ${\cR e}\,s > \theta$.\goodbreak

\section{Boundedness and Compactness of Composition Operators
$C_\phi \colon \cA^+ \to \cA^+$} 

\subsection{General results}

We begin by sharpening Newman's result ((a)
of the Introduction), under the form of the following (where it is assumed that 
$\phi$ is non-constant):

\begin{proposition}
The composition operator $C_\phi \colon A^+ \to A^+$ is compact if and
only if $\|\phi\|_\infty = \sup\limits_{z \in \DD} |\phi(z)| < 1$.
\end{proposition}

\noindent{\bf Proof.} As will be apparent from the proof of the next Proposition,
$C_\phi \colon A^+ \to A^+$ is compact if and only if
$\|\phi^n\|_{A^+} \to 0$ as $n \to \infty$.  On the other hand, by the spectral 
radius formula, we have $\|\phi\|_\infty = \lim\limits_{n \to \infty}
\|\phi^n\|_{A^+}^{1/n} = \inf\limits_{n \geq 1} \|\phi^n\|_{A^+}^{1/n}$. That 
finishes the proof.\par  
Alternatively, we could have applied to $f_n(z)=z^n$ a general criterion of 
Shapiro \cite{S}: 
``\emph{$C_\phi$ is compact if and only if 
$\|C_\phi(f_n)\|_{A^+} \to 0$ for each sequence $(f_n)_n$ in
$A^+$ which is bounded in norm and converges uniformly to zero on
compact subsets of\/ $\DD$}''.\hfill$\square$
\medskip

We now turn to the study of composition operators 
$C_\phi \colon \cA^+\to \cA^+$ associated with an analytic function
$\phi \colon \CC_0 \to \CC_0$.\par
We first recall the following:

\begin{theoreme}
\label{theoremnouveau}
(\cite[Theorem 4]{GH}).  Let $\phi \colon \CC_0 \to \CC$ be an analytic function 
such that $k^{-\phi} \in \cD$ for $k=1,2,\ldots\,$. Then we have necessarily:
\begin{equation}
\phi(s) = c_0s + \varphi(s), \mbox{~with~} c_0 \in \NN_0
\mbox{~and~} \varphi \in \cD. \label{equation1}
\end{equation}
\end{theoreme}

We will therefore restrict ourselves, in the sequel, to symbols
$\phi$ of the form given by (\ref{equation1}).  To avoid
trivialities, we will also assume once and for all that $\phi$ is
non-constant.
\goodbreak

\begin{theoreme}
\label{prop22}
Let $\phi \colon \CC_0 \to \CC$ be  an analytic function of the form
$(\ref{equation1})$.  Then~:
\begin{enumerate}
\renewcommand{\labelenumi}{(\alph{enumi})} 
\item
\begin{description}
\item {(i)} if $C_\phi$ maps $\cA^+$ into itself then 
$n^{-\phi} \in \cA^+$ and  $\|n^{-\phi}\|_{\cA^+} \leq C$, $n = 1,2,\ldots$, for a positive constant $C$ 
independent of $n$; 
\item {(ii)} conversely, if $(n^{-\phi})^\infty_{n=1}$ is a bounded sequence in $\cA^+$, then 
$\phi$ maps $\CC_0$ into $\CC_0$
and $C_\phi$ is a bounded composition operator on $\cA^+$.
\end{description} 
\renewcommand{\labelenumi}{(\alph{enumi})}
\item
\begin{description}
\item {(i)} $C_\phi \colon \cA^+ \to \cA^+$ is compact if and only if 
$\|n^{- \phi}\|_{\cA^+}\mathop{\longrightarrow}\limits_{n\to\infty} 0$. Then 
$\phi(\CC_0)\subseteq \CC_\delta$ for some $\delta>0$.
\item {(ii)} Assume that $\phi(s)=c_0s+\sum_1^\infty c_n n^{-s}$, with 
$\sum_1^\infty |c_n|<+\infty$. Then $C_\phi$ is compact if and only if 
$\phi(\CC_0)\subseteq \CC_\delta$ for some $\delta>0$.
\end{description} 
\end{enumerate}
\end{theoreme}

\noindent{\bf Proof.}
({\it a}) (\textit{i}) Suppose that $C_\phi$ maps $\cA^+$ into itself. $C_\phi$ is an
algebra homomorphism and $\cA^+$ is semi-simple, therefore (see \cite[p. 263]{R})
$C_\phi$ is continuous. Thus:
\begin{displaymath}
\|n^{-\phi}\|_{\cA^+} = \|C_\phi(n^{-s})\|_{\cA^+} \leq 
\|C_\phi\|\, \|n^{-s}\|_{\cA^+} = \|C_\phi\| = :C.
\end{displaymath}

(\textit{ii}) Conversely, suppose that $n^{-\phi} \in \cA^+$ and 
$\|n^{-\phi}\|_{\cA^+} \leq C$, 
$n = 1,2,\ldots$. We first see that, for $s \in \CC_0$, we have:
$n^{-{\cR e}\, \phi(s)} = |n^{-\phi(s)}| \leq \|n^{-\phi}\|_\infty \leq
\|n^{-\phi}\|_{\cA^+} \leq C$, whence 
${\cR e}\, \phi(s) \geq - \frac{\log C}{\log n}\cdot$ Letting $n$ tend to infinity gives
${\cR e}\, \phi(s) \geq 0$, and the open mapping theorem gives ${\cR e}\, \phi(s)> 0$, 
since $\phi$ is not constant.  If now 
$f(s) =\sum_1^\infty a_n n^{-s} \in \cA^+$, the series
$\sum_1^\infty a_n n^{- \phi (s)}$ is absolutely convergent
in $\cA^+$, so that $f \circ \phi \in \cA^+$, with
$\|f \circ \phi\|_{\cA^+} \leq \sum_1^\infty |a_n| \|n^{- \phi}\|_{\cA^+}
\leq C \sum_1^\infty |a_n| = C\|f\|_{\cA^+}$.\par\medskip 

(\textit{b}) (\textit{i}) Suppose that $C_\phi \colon \cA^+ \to \cA^+$ is compact.  Let
$f \in \cA^+$ be a cluster point of $n^{-\phi(s)} = C_\phi(n^{-s})$,
and let $(n_k)_k$ be a sequence of integers such that 
$\|n_k^{-\phi} - f\|_{\cA^+}\to 0$. For fixed $s \in \CC_0$, we have
$| n_k^{-\phi(s)} - f(s) | \leq \|n_k^{-\phi} - f\|_{\cA^+}$. But
$n_k^{-\phi(s)} \to 0$ (since ${\cR e}\, \phi(s)>0$, by part (a)), so that $f(s) = 0$. 
Hence $f = 0$. This implies $\|n^{-\phi}\|_{\cA^+} \to 0$.\par 
Now, since $\|n^{-\phi}\|_\infty\leq \|n^{-\phi}\|_{\cA^+}$, we get 
$n^{- \inf{\cR e}\,\phi(s)}=\|n^{-\phi}\|_\infty \to 0$, and so 
$\inf\limits_{s \in \CC_0}{\cR e}\,\phi(s) > 0$.\par  
Conversely, suppose that 
$\epsilon_n =\|n^{-\phi}\|_{\cA^+} \to 0$ and set 
$\delta_n = \sup_{k > n}\epsilon_k$. Let $T_n \colon \cA^+ \to \cA^+$ be the finite-rank
operator defined by $(T_n f)(s) = \sum_{k=1}^n a_k k^{-\phi(s)}$ if
$f(s) = \sum_{k=1}^\infty a_k k^{-s}$.  We have:
\begin{displaymath}
\|C_\phi f - T_nf \|_{\cA^+} \leq \sum\limits_{k > n} |a_k| \|k^{-\phi}\|_{\cA^+}
\leq \delta_n \sum\limits_{k > n} |a_k| \leq \delta_n \|f\|_{\cA^+},
\end{displaymath}
showing that $\|C_\phi - T_n\| \leq \delta_n$, and therefore that $C_\phi$
is compact.\par 

(\textit{ii}) For any $\upsilon\in {\cA^+}$, and for any real number $r\geq 1$, we have: 
\begin{equation}\label{etoile}
\| r^{-\upsilon}\|_{\cA^+}\leq r^{\|\upsilon\|_{\cA^+}}\,.
\end{equation}
Indeed:
\begin{displaymath}
r^{-\upsilon}=\exp(-\upsilon \log r)=\sum_{k=0}^\infty 
\frac{(-\log r)^k}{k!}\,\upsilon^k\in \cA^+
\end{displaymath}
since $\upsilon$ belongs to the algebra $\cA^+$. Moreover:
\begin{displaymath}
\|r^{-\upsilon}\|_{\cA^+}\leq \sum_{k=0}^\infty 
\frac{(\log r)^k}{k!}\,\|\upsilon\|_{\cA^+}^k 
=r^{\|\upsilon\|_{\cA^+}}
\end{displaymath}
(we may remark that when $\upsilon(s)=c_j j^{-s}$ is a monomial, we have equality; 
in particular: 
$\|n^{-c_j j^{-s}}\|_{\cA^+}=n^{|c_j|}$ for every positive integer $n$).
\medskip

We shall use the following:
\goodbreak

\begin{proposition}\label{Fait}
(see \cite{GH}) 
Let $\theta$ and $\tau$ be real numbers and suppose that $\phi$ maps $\CC_\theta$ into 
$\CC_\tau$. Then, if $\phi(s)=c_0 s +\varphi(s)$, and $\varphi$ is not constant, 
$\varphi$ maps $\CC_\theta$ into $\CC_{\tau-c_0 \theta}$.
\end{proposition}

Now, assume that $\varphi$ is non-constant (since otherwise the result is trivial), and 
that $\varepsilon=\inf\limits_{s\in \CC_0} {\cR e}\,\phi(s)>0$. By Proposition \ref{Fait}, 
$\varphi$ maps $\CC_0$ into $\CC_\varepsilon$. The spectral radius formula and Bohr's theory 
(as seen in the Introduction) give, with 
$\psi=2^{-\varphi}$:
\begin{displaymath}
\lim_{j\to+\infty} \|\psi\,^j\|_{\cA^+}^{1/j}=\sup_{h\in {\rm sp}\,\cA^+} | h(\psi)| 
= \sup_{s\in \CC_0} |\psi(s)|=2^{-\varepsilon};
\end{displaymath}
and, in particular, 
$\| 2^{-j\varphi}\|_{\cA^+}\mathop{\longrightarrow}\limits_{j\to+\infty} 0$.
Now, if $n$ is any positive integer, let $j=j(n)$ be the integer such that 
$2^j\leq n<2^{j+1}$, and set $r=n\,2^{-j}$, so that $1\leq r<2$. By using (\ref{etoile}), 
we get: 
\begin{align*}
\| n^{-\phi}\|_{\cA^+}&=\|n^{-\varphi}\|_{\cA^+}=\|2^{-j\varphi} r^{-\varphi}\|_{\cA^+}
\leq \|2^{-j\varphi}\|_{\cA^+}\,\|r^{-\varphi}\|_{\cA^+} \cr
&\leq 
\|2^{-j\varphi}\|_{\cA^+}\,r^{\|\varphi\|_{\cA^+}} 
\leq \|2^{-j\varphi}\|_{\cA^+}\,2^{\|\varphi\|_{\cA^+}}.
\end{align*}
This shows that $\|n^{-\phi}\|_{\cA^+}\mathop{\longrightarrow}\limits_{n\to\infty}0$ (more 
precisely, we have $\|n^{-\phi}\|_{\cA^+}={\rm O}\,(n^{-\delta})$ for some $\delta>0$), and so 
$C_\phi$ is compact, by part \textit{(b) (i)} of the theorem.\hfill $\square$
\bigskip 

\noindent{\bf Remark.} Using the notation of Theorem \ref{theoremnouveau}, we have:
\begin{displaymath}
\| n^{-\phi}\|_{\cA^+}=\| n^{-\varphi}\|_{\cA^+},
\end{displaymath}
and, in particular, the integer $c_0$ plays no role for the continuity or the compactness 
of the composition operator $C_\phi$ on $\cA^+$. This is quite amazing, since $c_0$ intervenes 
decisively in the study of composition operators on the Hilbert space $\cH^2$ of the 
square-summable Dirichlet series (so much so that Gordon and Hedenmalm \cite{GH} called it 
``\emph{characteristic}'').
\medskip

\begin{corollaire}
\label{prop23}
Let $\phi(s) = c_0 s + \sum\limits_{n=1}^\infty c_n n^{-s}$. Then $C_\phi$ is bounded
if ${\cR e}\,c_1 \geq \sum\limits_{n=2}^\infty |c_n|$, and is compact if 
${\cR e}\,c_1 > \sum\limits_{n=2}^\infty |c_n|$.
\end{corollaire}

\noindent{\bf Proof.} Let $\varphi_0\in \cA^+$ be defined by 
$\varphi_0(s)=\sum_{n=2}^\infty c_n n^{-s}$. 
For each positive integer $N$, we have:
$N^{-\phi(s)}=(N^{c_0})^{-s} N^{-c_1} N^{-\varphi_0(s)}$, and so 
the inequality (\ref{etoile}) with $r=N$ gives:
\begin{displaymath}
\|N^{-\phi}\|_{\cA^+}=N^{-{\cR e}\,c_1}\|N^{-\varphi_0}\|_{\cA^+} 
\leq N^{-{\cR e}\,c_1} N^{\|\varphi_0\|_{\cA^+}} 
=N^{-{\cR e}\,c_1+\sum_{n=2}^\infty |c_n|};
\end{displaymath}
thus $\|N^{-\phi}\|_{\cA^+}$ is less than $1$ in the first case, and tends to $0$ in the 
second case. Theorem \ref{prop22} ends the proof.\hfill$\square$ 
\medskip  

Note that under the assumption of Corollary \ref{prop23}, $C_\phi\colon \cA^+\to \cA^+$ is 
actually a contraction: $\|C_\phi\| \leq 1$. 
\medskip 

\subsection{Some specific examples}

One of the main differences between the study of composition operators on 
$\cA^+$ and those on $A^+(\TT)$ is the fact that the function $z \mapsto z$ does not 
belong to $\cA^+$. Therefore, it is not clear that if $C_\phi$ is a composition 
operator on $\cA^+$, we must have $\sum_n | c_n|<+\infty$. In some cases, 
it is however true. The next proposition contains a partial result of this type.
\medskip\goodbreak

Recall (see \cite{K}) that $(\lambda_j)_{j\geq 1}$ is a Sidon set if
\begin{displaymath}
\sum_{j=1}^N |a_j| \leq C_0 \sup_{t \in \RR} 
\Big|\sum_{j=1}^N a_j e^{i \lambda_j t} \Big|
\end{displaymath}
for some finite positive constant $C_0$. 

\begin{proposition}\label{prop4bis} $\ $
\begin{enumerate}
\renewcommand{\labelenumi}{(\alph{enumi})} 
\item If $2 \leq q_1 < q_2 < \cdots$ are multiplicatively independent integers and 
$\phi(s) = c_0 s + c_1 +\sum_{j=1}^\infty d_j q^{-s}_j$, then the boundedness of 
$C_\phi\colon \cA^+\to \cA^+$ implies  
${\cR e}\,c_1 \geq \sum_{j=2}^\infty |d_j|$, and its compactness implies 
${\cR e}\,c_1 > \sum_{j=2}^\infty |d_j|$.\par
\item  Let $(\lambda_j)_{j \geq 1}$ be a Sidon set of positive
integers, $r$ an integer $\geq 2$, and 
$\phi(s) = c_0 s +\varphi(s)$, where $\varphi \in {\cD}$ and $\varphi(s) = c_1 +
\sum_{j=1}^\infty d_j r^{-\lambda_j s}$ for ${\cR e}\,s $ large. 
Then the boundedness of $C_\phi\colon {\cA}^+\to {\cA}^+$ requires that 
$\sum_{j=1}^\infty |d_j|<+\infty$.
\end{enumerate}
\end{proposition}

\noindent{\bf Proof.}
\textit{(a)} Write $\varphi_0(s)=\sum_{j=1}^\infty d_j q_j^{-s}$, as in the proof of 
Corollary \ref{prop23}. For every integer $n\geq 2$, we have, for ${\cR e}\,s$ large enough:
\begin{displaymath}
n^{-\phi(s)}=(n^{c_0})^{-s} n^{-c_1}\exp\big(-\varphi_0(s)\,\log n\big)
=(n^{c_0})^{-s} n^{-c_1}\sum_{k=0}^\infty \frac{(-\log n)^k}{k!}\,\varphi_0(s)^k.
\end{displaymath}
Since $C_\phi$ is assumed to be bounded on $\cA^+$, we know that $n^{-\phi}\in \cA^+$, and 
so:
\begin{displaymath}
n^{-\phi(s)}=\sum_{j=1}^\infty a_{n,j} j^{-s},\hskip 3mm
\mbox{with}\ \sum_{j=1}^\infty |a_{n,j}|<+\infty.
\end{displaymath}
But the supports (spectra) of the $\varphi_0^k$'s do not intersect: in fact, the spectrum of 
$\varphi_0^k$ only involves finite products $\prod_j q_j^{\alpha_j}$, where $\sum \alpha_j=k$, 
and these products are all distinct. In particular, for $k=1$, 
$(-\log n)\varphi_0(s)$ is part of the expansion of $n^{-\phi(s)}$, 
which means that $(-d_j\log n)_j$ is a subsequence of $(a_{n,j})_j$. Therefore, 
$\sum_j |d_j|<+\infty$ (and so $\varphi_0\in \cA^+$), and the series expansion of 
$\varphi_0(s)$ holds for every $s\in \CC_0$.\par\noindent
Finally, since the $\log q_j$'s are rationally independent, Kronecker's Approximation Theorem 
implies that, for each $\sigma>0$, we have:
\begin{displaymath}
\inf_{t\in \RR} {\cR e}\, \phi(\sigma+it)=c_0\sigma +{\cR e}\, c_1 -
\sum_{j=1}^\infty |d_j| q_j^{-\sigma}.
\end{displaymath}
Since the left-hand side is $\geq 0$ by the first part of Theorem \ref{prop22}, we get 
${\cR e}\,c_1\geq \sum_{j=1}^\infty |d_j|$, by letting $\sigma$ go to zero. The 
compact case is similar.
\medskip

\textit{(b)} We have (see \cite{LR}):
\begin{equation}\label{S1}
\inf_{\tau \in \RR} \sum_{j=1}^N \rho_j 
\cos(\lambda_j \tau + \xi_j) \leq - \delta
\sum_{j=1}^N \rho_j
\end{equation}
for some other constant $\delta > 0$, where the $\rho_j$'s
(non-negative) and the (real) $\xi_j$'s are arbitrary.
Without loss of generality, we can assume that $r=2$.  Fix an
integer $J \geq 1$, and let $B : \RR \to \RR^+$ (see
\cite{Kat}, p.165) be a non-negative Dirichlet polynomial
(of the form $\sum \alpha_k e^{i \beta_k t}$, $\beta_k \in \RR$,
$\alpha_k \in \CC)$ such that:
\begin{equation}\label{S2}
\widehat B(0) = \widehat B(\lambda_j \log 2) = 1,\hskip 3mm 1 \leq j \leq J
\end{equation}
(recall that $\widehat B(\lambda) = \lim\limits_{T \to \infty}
\frac{1}{2T} \int^T_{-T} B(t) e^{-i \lambda t}\,dt$). \par
For large $\sigma>0$, we have an absolutely convergent expansion:
\begin{displaymath}
\varphi(\sigma + i(t + \tau)) = c_1 + \sum_{j=1}^\infty
d_j 2^{-\lambda_j \sigma} 2^{-\lambda_j it} e^{-i\lambda_j \tau\log 2},
\end{displaymath}
so that, for ${\cR e}\,s$ large enough (say ${\cR e}\,s\geq \sigma_0>0$):
\begin{equation}\label{S3}
\lim_{T \to \infty} \frac{1}{2T} 
\int^T_{-T}\!\!\varphi (s + i\tau) B(\tau)\,d\tau =
c_1 + \sum_{j=1}^\infty d_j 2^{-\lambda_j s} \widehat B(\lambda_j \log 2).
\end{equation}
Actually, (\ref{S3}) holds for every $s$ with positive real part $\sigma$. To see this, 
set: 
\begin{displaymath}
f_T(s)=\frac{1}{2T}\int_{-T}^T \varphi (s+i\tau) B(\tau)\,d\tau.
\end{displaymath}
Proposition \ref{Fait} shows that ${\cR e}\,\varphi(s+i\tau)>0$ for $s\in \CC_0$, 
and thus that ${\cR e}\, f_T(s)>0$ for $s\in \CC_0$. Moreover, $f_T$ as well as the 
right-hand side of (\ref{S3}), since $B$ is a Dirichlet polynomial, are holomorphic 
in $\CC_0$; hence a normal family argument gives the above statement.\par  
Therefore, if we take the real part of both sides of (\ref{S3}), we get, for every 
$\sigma>0$ and $\t\in\RR$:
\begin{displaymath}
{\cR e}\,c_1  + \sum_{j\geq 1} 2^{-\lambda_j \sigma} {\cR e}\,
\big(d_j 2^{-\lambda_j it}  \widehat B  (\lambda_j \log 2)\big) 
= \lim_{T \to +\infty} {\cR e}\,f_T(\sigma+it)\geq 0.
\end{displaymath}
Letting $\sigma$ tend to zero gives:
\begin{displaymath}
{\cR e}\,c_1 + \sum_{j \geq 1} {\cR e}\, \big(d_j 2^{-\lambda_j it} 
\widehat B (\lambda_j \log 2)\big) \geq 0, \quad
\mbox{for any $t \in \RR.$}
\end{displaymath}
Taking the infimum with respect to $t$ and using (\ref{S1}), we get:
\begin{displaymath}
{\cR e}\,c_1 - \delta \sum_{j=1}^\infty |d_j|~|\widehat B (\lambda_j \log 2)| \geq 0
\end{displaymath}
and therefore, using (\ref{S2}):
\begin{displaymath}
{\cR e}\,c_1 - \delta \sum_{j=1}^J |d_j| \geq 0.
\end{displaymath}
\indent It follows that $\sum_{j=1}^\infty |d_j| \leq \frac{1}{\delta} {\cR e}\,c_1$, and this 
ends the proof of Proposition \ref{prop4bis}.\hfill$\square$
\bigskip

\noindent{\bf Remark.} The above proof gives the following information about 
Dirichlet series, which is actually not connected to composition operators: \emph{let 
$\varphi$ be a Dirichlet series which can be written as 
$\varphi(s)=c_1+\sum_{j\geq 1} d_j r^{-\lambda_j s}$, where $(\lambda_j)_j$ is a Sidon 
sequence; if there is a $\beta\in\RR$ such that $\varphi(\CC_0)\subseteq \CC_\beta$, then 
$\sum_{j\geq 1} |d_j|<+\infty$.}
\bigskip

However, in general, conditions like $\sum_{n\geq 2} |c_n|\leq {\cR e}\,c_1$ 
(resp. $<{\cR e}\,c_1$) are not necessary to have boundedness or compactness of the 
composition operator $C_\phi\colon \cA^+\to \cA^+$ (with 
$\phi(s)=c_0 s+ c_1+\sum_{n\geq 2} c_n n^{-s}$), 
as shown by the following examples.

\begin{proposition}
\label{prop24}
Let $\phi(s) = c_0 s + c_1 + c_r r^{-s} + c_{r^2} r^{-2s}$, where $r \geq 2$
and $c_r$, $c_{r^2}$ are $> 0$. Then:
\begin{enumerate}
\renewcommand{\labelenumi}{(\alph{enumi})} 
\item If we have
\begin{eqnarray}
{\cR e}\, c_1 > \frac{(c_r)^2}{8 c_{r^2}} + c_{r^2}, \label{equation4}
\end{eqnarray}
$C_\phi \colon \cA^+ \to \cA^+$ is bounded and even compact.
\item Conversely, if $C_\phi \colon \cA^+ \to \cA^+$ is bounded, and
moreover $c_r \leq 4c_{r^2}$, we must have
\begin{eqnarray}
{\cR e}\, c_1 \geq \frac{(c_r)^2}{8 c_{r^2}} + c_{r^2}.
\label{equation4'}
\end{eqnarray}
In fact, we must have (\ref{equation4}) whenever $C_\phi$ is compact.
\item If
\begin{eqnarray}
{\cR e}\, c_1 = \frac{(c_r)^2}{8 c_{r^2}} + c_{r^2},
\label{equation4''}
\end{eqnarray}
then $C_\phi\colon \cA^+\to \cA^+$ is bounded if and only if $c_r\neq 4c_{r^2}$. 
\end{enumerate}
\end{proposition}

\noindent{\bf Proof.} \textit{(a)} and \textit{(b)} follow immediately from 
Theorem \ref{prop22}, since (\ref{equation4}) 
implies ${\cR e}\,\phi(s)> c_0 {\cR e}\,s +\delta\geq \delta$ for every $s\in \CC_0$, with 
$\delta={\cR e}\, c_1 - \big[\frac{(c_r)^2}{8 c_{r^2}} + c_{r^2}\big]$, and, under the 
assumption that $c_r \leq 4c_{r^2}$, the converse is true.\hfill$\square$
\medskip

However, we shall give another proof, because we think that it sheds additional light.
\medskip

\noindent{\bf Second proof.}\par
\textit{(a)} Without loss of generality, we may and shall assume that $r=2$.  We will make use 
(see \cite[p. 60]{L}) of the Hermite polynomials $H_0,H_1,\ldots$ defined by:
\begin{eqnarray}
H_k(\lambda) = (-1)^k e^{\lambda^2} \frac{d^k}{d\lambda^k} \left(e^{-\lambda^2}
\right) =
(2\lambda)^k + \mbox{terms of lower degree.} \label{equation5}
\end{eqnarray}
The exponential generating function of the $H_k$'s is:
\begin{eqnarray}
\sum^\infty_{k=0} \frac{H_k(\lambda)}{k!} x^k =
\exp \left(2 \lambda x - x^2 \right). \label{equation6}
\end{eqnarray}
Following Indritz \cite{I}, we have the sharp estimate
\begin{eqnarray}
|H_k(\lambda)| \leq (2^k k!)^{1/2} e^{\lambda^2/2}
\label{equation7},
\end{eqnarray}
for each $k \in \NN_0$ and each $\lambda \in \RR$.
The estimate (\ref{equation7}) implies the following:

\begin{lemme} \label{lem5}
Let $\lambda$ be a real number, and $x$ be a
non-negative real number.  Then we have:
\begin{eqnarray}
\sum^\infty_{k=0} \frac{|H_k(\lambda)|}{k!} x^k \leq C(1+x)^{1/2}
\exp \left(x^2 + \frac{\lambda^2}{2} \right)
\label{equation8}
\end{eqnarray}
where $C$ is a positive constant.
\end{lemme}

\noindent{\bf Proof of the Lemma.} (\ref{equation7}) implies that:
\begin{displaymath}
\sum_{k=0}^\infty \frac{|H_k(\lambda)|}{k!} x^k \leq
\sum_{k=0}^\infty \frac{(x \sqrt 2)^k}{(k!)^{1/2}} e^{\lambda^2/2}.
\end{displaymath}
We now make use of the classical estimate (see \textit{e.g.} Dieudonn\'e
\cite[p. 195]{D}):
\begin{eqnarray}
\sum_{k=0}^\infty \frac{y^k}{(k!)^p} \sim \frac{1}{\sqrt p}
(2\pi)^{\frac{1-p}{2}} y^{\frac{1-p}{2p}} \exp(py^{1/p})
\mbox{~~as $y \to \infty$ ($p > 0$ fixed).} \label{equation9}
\end{eqnarray}
Using the above, with $p = \frac 12$
and $y = x \sqrt 2$, we obtain for some constant~$C$~:
\begin{displaymath}
\sum_{k=0}^\infty  \frac{|H_k(\lambda)|}{k!} x^k \leq C
e^{\lambda^2/2}(1+x)^{1/2} e^{x^2},
\end{displaymath}
proving the Lemma.\hfill$\square$
\bigskip

Note that, if we wish to avoid the use of (\ref{equation9}), we can easily 
obtain the slightly weaker estimate:
\begin{eqnarray}
\sum_{k=0}^\infty \frac{|H_k(\lambda)|}{k!} x^k \leq C_a \exp
\left(ax^2 + \frac{\lambda^2}{2} \right),\hskip 3mm \mbox{~for each $a >1$.}
\label{equation8'}
\end{eqnarray}
Indeed, we have by the Cauchy-Schwarz inequality and by (\ref{equation7})~:
\begin{eqnarray*}
\sum_{k=0}^\infty \frac{|H_k(\lambda)|}{k!} x^k
 & = & \sum_{k=0}^\infty \frac{|H_k(\lambda)|}{(k!)^{1/2} (2a)^{k/2}}
  \frac{(2a)^{k/2} x^k}{(k!)^{1/2}} \\
 & \leq & \left( \sum_{k=0}^\infty
 \frac{|H_k(\lambda)|^2}{k!(2a)^k} \right)^{1/2}
 \left( \sum_{k=0}^\infty \frac{ (2a)^k x^{2k}}{k!} \right)^{1/2} \\
 & \leq & e^{\lambda^2/2} \Big( \sum_{k=0}^\infty a^{-k}\Big)^{1/2} \exp (ax^2) \\
 & = & (1-a^{-1})^{-1/2} \exp \left(ax^2 + \frac{\lambda^2}{2}\right).
\end{eqnarray*}

We now finish the proof of Proposition \ref{prop24}. First, we notice that:
\begin{displaymath}
n^{-\phi(s)} = (n^{c_0})^{-s} n^{-c_1} \exp \left(-c_2 2^{-s} \log n
- c_4 4^{-s} \log n\right).
\end{displaymath}  
We then set:
\begin{eqnarray}
x_n = \sqrt{c_4 \log n},\hskip 3mm \lambda_n = \frac{-c_2}{2 \sqrt{c_4}}
\sqrt{\log n},\hskip 3mm x = 2^{-s} x_n, \label{equation10}
\end{eqnarray}
which allows us to write $n^{-\phi(s)}$ under the form:
\begin{displaymath}
n^{-\phi(s)} = (n^{c_0})^{-s} n^{-c_1} \exp
\big(2 \lambda_n x - x^2 \big) 
= (n^{c_0})^{-s} n^{-c_1}
\sum_{k=0}^\infty \frac{H_k(\lambda_n)}{k!} x_n^k(2^k)^{-s}.
\end{displaymath}
This implies that we have the equality:
\begin{eqnarray}
\|n^{-\phi}\|_{\cA^+} = n^{-{\cR e}\, c_1} \sum_{k=0}^\infty
\frac{|H_k(\lambda_n)|}{k!} x_n^k \label{equation11}.
\end{eqnarray}
If we now use Lemma \ref{lem5} and change $C$ (if necessary), we get
for $n \geq 2$~:
\begin{align*}
\|n^{-\phi}\|_{\cA^+} &\leq  C n^{-{\cR e}\, c_1} (\log n)^{1/4} \exp
\big(x^2_n + \frac{\lambda^2_n}{2} \big) \cr
& = C(\log n)^{1/4} n^{-{\cR e}\, c_1} n^{ \frac{c_2^2}{8c_4} + c_4} 
= : \epsilon_n.
\end{align*}
By (\ref{equation4}), we have $\epsilon_n \to 0$, implying that
$C_\phi \colon \cA^+ \to \cA^+$ is compact as a consequence of Theorem 
\ref{prop22}.\smallskip

\textit{(b)} The identities (\ref{equation11}) and (\ref{equation6})
imply that we have, for each real $\theta$,
\begin{eqnarray*}
n^{{\cR e}\, c_1} \|n^{-\phi}\|_{\cA^+} & \geq & \left| \sum_{k=0}^\infty
\frac{H_k(\lambda_n)}{k!} x_n^k e^{ik \theta} \right| =
\left| \exp \left(2 \lambda_n x_n e^{i \theta} - x_n^2
e^{2i\theta}\right)\right| \\
& = & \exp \left(2 \lambda_n x_n \cos \theta - x_n^2 \cos 2 \theta
\right).
\end{eqnarray*}
Setting \ $t = \cos \theta$, we see that 
\begin{displaymath}
2 \lambda_n x_n \cos \theta
- x_n^2 \cos 2 \theta = 2 \lambda_n x_n t - x_n^2(2t^2 - 1)
\end{displaymath}
is maximum for $t = \frac{\lambda_n}{2x_n} = \frac{-c_2}{4c_4}$, and
this $t$ will be admissible if $\left| \frac{-c_2}{4c_4}\right| \leq
1$, \textit{i.e.} if $c_2 \leq 4c_4$ (recall that $c_2,c_4$ are positive).
For this value of $t$, we get:
\begin{eqnarray}
\|n^{-\phi}\|_{\cA^+} \geq n^{-{\cR e}\, c_1 + \frac{c^2_2}{8c_4} + c_4},\hskip 3mm 
n = 1,2,\ldots \hskip 3mm,\quad c_2 \leq 4c_4. \label{equation11prime}
\end{eqnarray}
Now, if $C_\phi$ is bounded, $\|n^{-\phi}\|_{\cA^+}$ is bounded from
above, and (\ref{equation11prime}) implies that
${\cR e}\, c_1 \geq \frac{(c_2)^2}{8c_4} + c_4$.  If $C_\phi$ is
compact, $\|n^{-\phi}\|_{\cA^+} \to 0$ and (\ref{equation11prime}) implies
that ${\cR e}\, c_1 > \frac{(c_2)^2}{8c_4} + c_4$.\hfill $\square$
\bigskip

\noindent{\bf Remark.} Condition (\ref{equation4}) is 
a more general sufficient condition for the boundedness of $C_\phi$
than the ``\emph{trivial}'' sufficient condition
${\cR e}\, c_1 \geq |c_2| + |c_4|$ of Corollary \ref{prop23} if and only 
if $c_r < 8 c_{r^2}$. This might be due to the highly
oscillatory character of the Hermite polynomials $H_k(\lambda)$,
involving a term $\cos \left(\sqrt {2k+1} \lambda - k \frac{\pi}{2} \right)$
(see \cite[p. 67]{L}), which we ignore when we majorize $|H_k(\lambda)|$
as in (\ref{equation7}).
\bigskip 

\noindent{\bf End of proof of Proposition \ref{prop24}.} \textit{(c)} We still assume 
that $r=2$. First, if $c_2>4 c_4$, then (\ref{equation4''}) implies that 
$\phi(\CC_0)\subseteq \CC_\delta$ for some $\delta>0$, and we are done. So, we 
assume that $c_2\leq 4 c_4$. We have:
\begin{align*}
\| (2^j)^{-\phi}\|_{\cA^+} 
& = \| 2^{-j\phi}\|_{\cA^+}= \| (2^{-c_1 -c_2 2^{-s} - c_4 4^{-s}})^j\|_{\cA^+} \\
&= \big\|\big(\exp[(-c_1 -c_2 2^{-s} -c_4 4^{-s})\log 2]\big)^j\big\|_{\cA^+} 
= \|\psi^j\|_{A^+(\TT)},
\end{align*}
with
\begin{displaymath}
\psi(z)=\exp\big(-(c_1+c_2 z + c_4 z^2)\log 2\big).
\end{displaymath}
We then apply Newman's result (quoted as (a) in the Introduction: see \cite{Ne}) to 
check whether the sequence $(\|\psi^j\|_{A^+(\TT)})_j$ is bounded. Let 
$\theta_0\in [0,2\pi[$ be such that $|\psi(e^{i\theta_0})|=1$. We look for the 
coefficient of $t^2$ in the Taylor expansion of:
\begin{displaymath}
\log \psi(e^{i\theta_0+it})= 
-(c_1+c_2 e^{i\theta_0}e^{it}+c_4 e^{2i\theta_0}e^{2it})\log 2.
\end{displaymath}
This term is:
\begin{displaymath}
\Big(\frac{c_2}{2} e^{i\theta_0} + 2 c_4 e^{2i\theta_0}\Big)\log 2,
\end{displaymath}
and its real part is:
\begin{equation}\label{equation4T}
\Big(\frac{c_2}{2} \cos \theta_0 + 2 c_4 (2\cos^2\theta_0 -1)\Big)\log 2.
\end{equation}
Now, remark that the condition $|\psi(e^{i\theta_0})|=1$ means that:
\begin{displaymath}
{\cR e}\,c_1=-c_2 \cos \theta_0 - c_4 (2\cos^2\theta_0 -1),
\end{displaymath}
which gives, using (\ref{equation4''}), $(4 c_4 \cos\theta_0 +c_2)^2=0$, that is 
$\cos \theta_0= -\frac{c_2}{4 c_4}\cdot$ Hence (\ref{equation4T}) 
is equal to $0$ if and only if $c_2=4 c_4$.\par
But in this case, $\theta_0=\pi$, and Taylor's expansion becomes:
\begin{displaymath}
\log \psi(e^{i(\theta_0+t)})= d_1 +d_2 t + 0.t^2 + i\log 2 \frac{2c_4}{3}t^3+\cdots
\end{displaymath}
Hence, in Newman's terminology (see \cite{Ne}, and see (a) in the Introduction), the point 
$e^{i\theta_0}$ is not an ordinary point, and so the sequence 
$(\|\psi^j\|_{A^+(\TT)})_j$ is not bounded. It follows that 
the sequence $(\|2^{-j\phi}\|_{\cA^+})_j$ is not bounded either.\par
In the case $c_2<4 c_4$, the point $e^{i\theta_0}$ is ordinary, and so  
$(\|2^{-j\phi}\|_{\cA^+})_j$ is bounded. Since 
$\sum_n |c_n|=|c_1|+|c_2|+|c_4|<+\infty$, 
the argument used in the proof in Theorem \ref{prop22}, \textit{(b)~(ii)} gives the 
boundedness of $C_\phi$.\hfill$\square$
\medskip

\noindent{\bf Remark.} Part \textit{(c)} of Proposition \ref{prop24} shows that, if   
${\cR e}\,c_1=\frac{(c_r)^2}{8c_{r^2}} + c_{r^2}$ and $c_r=4c_{r^2}$ (so 
${\cR e}\,c_1=3 c_4$, and $\psi(s)=ia + c(3+4\cdot 2^{-s} + 4^{-s})$, with $a\in\RR$ 
and $c>0$), then $C_\phi$ is not bounded on $\cA^+$, though 
$\phi(\CC_0)\subseteq \CC_0$ (and $\sum_n |c_n|<+\infty$).\par

\section{Automorphisms of $A^+(\hbox{\large $\TT\,$}^k)$, 
$A^+(\hbox{\large $\TT\,$}^\infty)$, $\cA^+$} 

In this section, we will make repeated use of the following Lemma
(see (b) of the Introduction):

\begin{lemme}\label{lemme31}
Let $\phi(z) = \prod\limits^J_{j=1} \epsilon_j \frac{z-a_j}{1-\overline{a_j}z},$
where $|\epsilon_j| = 1$ and $a_j \in \DD$.  Suppose that $\|\phi^n\|_{A^+}$ 
remains bounded $(n = 1,2,\ldots)$.  Then,
$a_j = 0$ for each $j$.
\end{lemme}

\noindent{\bf Proof.} This Lemma is well-known (see \cite{Ne} page 38, assertion $(3)$,or \cite{K}, page 77).
For example, if $a_j \not=0 $ for some $j$, we have
$\phi(e^{it}) = e^{ig(t)}$, where $g$ is a $\cC^2$, real, non affine
function; and the Van der Corput inequalities show that we even
have: $\|\phi^n\|_{A^+} \geq \delta \sqrt n$.\hfill $\square$
\medskip

Since $|\phi(e^{it})| = 1$, Lemma \ref{lemme31} can
be viewed as a special case of the following Lemma (which will be
needed only in Section 4, but which we state here because it is
the natural extension of Lemma \ref{lemme31}), due to Beurling and
Helson, and this Lemma is itself a special case of Cohen's Theorem
(see \cite{R}, page 93, corollary of Theorem 4.7.3). We shall use the following definition: \par\smallskip
\noindent
Let $G$ be a discrete abelian group, and $\Gamma$ be its
(compact) dual group; the Wiener algebra $A(\Gamma)$ is the set
of functions $f : \Gamma \to \CC$ which can be written as an
absolutely convergent series $f(\gamma) = \sum_1^\infty a_n(x_n,\gamma)$, with the 
norm $\|f\|_{A(\Gamma)} = \sum_1^\infty |a_n|$, and where $(x_n,\gamma)$ 
denotes the action of 
$\gamma\in \Gamma$ on the element $x_n$ of $G$.  We are now ready to state:

\begin{lemme} \label{lemme8}
(Beurling-Helson).  Let $G$ be a discrete abelian group, with
\emph{connected} dual group $\Gamma$. Let $\phi \in A(\Gamma)$, which does 
not vanish on $\Gamma$, and such that $\|\phi^n\|_{A(\Gamma)} \leq C$ for
some constant $C$ $(n = 0,\pm 1,\pm 2,\ldots)$.  Then, $\phi$ is
affine, i.e. there exist a complex number $a$ with $|a| = 1$ and
an element $x$ of $G$ such that $\phi(\gamma) = a(x,\gamma)$ for
any $\gamma \in \Gamma$.
\end{lemme}

Let us now consider the Wiener algebra $A^+(\TT^k)$ in $k$
variables, \textit{i.e.} the algebra of functions $f : \overline \DD^k \to \CC$
which can be written as:
\begin{displaymath}
f(z) = \sum_{n_1,\ldots,n_k \geq 0} a(n_1,\ldots,n_k) z_1^{n_1}\ldots z_k^{n_k}, 
\hskip 3mm z = (z_1,\ldots,z_k),
\end{displaymath}
with the norm
$\|f\|_{A^+(\TT^k)} = \sum\limits_{n_1,\ldots,n_k \geq 0} |a(n_1,\ldots,n_k)| 
<+\infty$.\par  
If $\phi = (\phi_1, \ldots, \phi_k) : \DD^k \to \CC^k$ is an
analytic function, the composition operator $C_\phi$ will be
bounded on $A^+(\TT^k)$ if and only if :
\begin{eqnarray}
\|\phi^n_j\|_{A^+(\TT^k)} \leq C,\hskip 8mm j = 1, \ldots, k, \mbox{~and~}\hskip 3mm 
n = 0,1,2,\ldots
\label{equation31}
\end{eqnarray}
(the proof is the same as in Newman's case $k=1$).\par\smallskip

Then, since 
$\|\phi_j\|_\infty = \lim\limits_{n \to \infty} \|\phi_j^n\|_{A^+(\TT^k)}^{1/n}$, we see 
that $\phi$ necessarily maps $\DD^k$ into $\overline\DD^k$.  We can now state:

\begin{theoreme}
\label{theorem33} Assume that the map $\phi\colon \DD^k \to \overline \DD^k$ induces a
bounded operator \hbox{$C_\phi \colon A^+(\TT^k) \to A^+(\TT^k)$.}  Then $C_\phi$ is an 
automorphism of $A^+(\TT^k)$ if and only if 
$\phi(z) = (\epsilon_1 z_{\sigma(1)}, \ldots, \epsilon_k z_{\sigma(k)})$
for some permutation $\sigma$ of $\{1,\ldots,k\}$ and some complex
signs $\epsilon_1,\ldots,\epsilon_k$.
\end{theoreme}

\noindent{\bf Proof.} The sufficient condition is trivial.   
For the necessary one, we first observe that, for $j = 1, \ldots, k$, 
$\phi_j \in A^+(\TT^k)$, since $\phi_j = C_\phi(z_j)$; hence 
$\phi$ can be continuously extended to a continuous map,
still denoted by $\phi$, from $\overline \DD^k$ to $\overline \DD^k$.
We are going to show that this map is bijective.\\  
Assume first that $a,b\in \overline \DD^k$ and that $\phi(a) = \phi(b)$.  Let 
$f \in A^+(\TT^k)$; since $C_\phi$ is bijective we can find $g \in A^+(\TT^k)$
such that $f = g \circ \phi$, so that $f(a) = f(b)$.  Since
$A^+(\TT^k)$ obviously separates the points of $\overline \DD^k$, we
have $a=b$. In particular, $\phi$ is injective on $\DD^k$ and by Osgood's Theorem (see
\cite{Na}, Theorem 5, page 86) $\det(\phi'(z)) \not= 0$ for each $z \in \DD^k$, implying that $\phi$ is an 
open mapping on $\DD^k$.  Therefore, $\phi(\DD^k)\subseteq \DD^k$.\\ 
Now, let $u \in \overline \DD^k$.  Define an element $L$ of the spectrum of $A^+(\TT^k)$ 
by $L(f) = g(u)$ if $f = g \circ \phi$. Since the spectrum of $A^+(\TT^k)$ is clearly 
$\overline \DD^k$, we can find $a \in \overline \DD^k$ such that $L(f) = f(a)$, so that 
$g\big(\phi(a)\big) = g(u)$ for any $g \in A^+(\TT^k)$, implying $u = \phi(a)$. 
$\phi$ is therefore a homeomorphism~: $\overline \DD^k \to \overline \DD^k$.\par  
Since $\phi(\overline \DD^k) = \overline \DD^k$ and $\phi(\DD^k)\subseteq \DD^k$, we get 
$\phi(\DD^k) = \DD^k$. In particular, $\phi\in \mbox{Aut~} \DD^k$, the group of analytic 
automorphisms of $\DD^k$.\par \smallskip
 
Recall that (\cite{Na}, Proposition 3, page 68):

\begin{lemme}
\label{lemme33}
The analytic map $\phi\colon \DD^k \to \DD^k$ belongs to \mbox{\rm Aut~}$\DD^k$ if and
only if
$$\phi(z) = \left( \epsilon_1 \frac{z_{\sigma(1)} - a_1}{1 -
\overline{a_1} z_{\sigma(1)}}\raise0,5mm\hbox{,} \cdots\raise0,5mm\hbox{,} \epsilon_k
\frac{z_{\sigma(k)} - a_k}{1 - \overline{a_k} z_{\sigma(k)}}\right),$$
for some permutation $\sigma$ of $\{1,\ldots,k\}$, for some
$(a_1,\ldots,a_k) \in \DD^k$ and some complex signs $\epsilon_1,\ldots,
\epsilon_k$.
\end{lemme}

We therefore see that $\phi_j(z) = \epsilon_j
\frac{z_{\sigma(j)} - a_j}{1 - \overline{a_j} z_{\sigma(j)}}$, so
that for each $n \in \NN$, we have in view of (\ref{equation31}):
\begin{displaymath}
\left\| \left( \epsilon_j \frac{z - a_j}{1 - \overline{a_j} z}\right)^n \right\|_{A^+} = 
\left\| \left( \epsilon_j\frac{z_{\sigma(j)} - a_j}{1 - \overline{a_j} z_{\sigma(j)}} 
\right)^n\right\|_{A^+(\TT^k)} \leq C.
\end{displaymath}

Lemma \ref{lemme31} now implies that $a_j = 0$, $j = 1,\ldots, k$,
so that $\phi_j(z) = \epsilon_j z_{\sigma(j)}$, and this ends the
Proof of Theorem \ref{theorem33}.\hfill$\square$
\bigskip

We now consider the Wiener algebra $A^+(\TT^\infty)$
in countably many variables. It will be convenient to consider
holomorphic functions on the open unit ball 
${\textbf B} = \DD^\infty \cap {\textbf c}_0$
of the Banach space ${\textbf c}_0$ of sequences $z = (z_n)_{n \geq 1}$
tending to zero at infinity, with its natural norm 
$\|z\| =\sup_{n \geq 1} |z_n|$. We then have the following extension of
Cartan's Lemma \ref{lemme33} to the case of ${\textbf B}$, which is due to Harris
(\cite{Ha}):

\begin{lemme}\label{lemme34} {\bf (Analytic Banach-Stone Theorem)}.  The 
analytic automorphisms $\phi \colon {\textbf B} \to {\textbf B}$ are exactly the maps of 
the form $\phi = (\phi_j)_{j \geq 1}$, with 
$\phi_j(z) = \epsilon_j 
\frac{z_{\sigma(j)} - a_j}{1 - \overline{a_j} z_{\sigma(j)}}\,$\raise0,5mm\hbox{,} for 
some permutation $\sigma$ of $\NN$, some point 
$a = (a_j)_{j \geq 1}\in {\textbf B}$, 
and some sequence $(\epsilon_j)_{j\geq1}$ of complex signs.
\end{lemme}

Recall that the linear Banach-Stone Theorem states : if 
$L \colon {\textbf c}_0 \to {\textbf c}_0$
is a surjective isometry fixing the origin, then $L$ has the form:
\begin{displaymath}
L(z_1,\ldots,z_n,\ldots) = (\epsilon_1 z_{\sigma(1)}, \ldots,
\epsilon_n z_{\sigma(n)}, \ldots).
\end{displaymath}

If we want to exploit Lemma \ref{lemme34} for describing the
composition automorphisms of $A^+(\TT^\infty)$, we have to make an
extra assumption, the reason for which is the following: if $C_\phi$ is an automorphism of $A^+(\TT^\infty)$,
then $\phi$ is an automorphism of $\overline \DD^\infty$, but there
is no reason, \textit{a priori}, why $\phi$ should be an automorphism of
${\textbf B}$.  

\begin{theoreme}
\label{theo32} Let $\phi = (\phi_j)_j \colon {\textbf B} \to {\textbf B}$ be an analytic map
such that $C_\phi$ maps $A^+(\TT^\infty)$ into itself.  Then~:
\begin{enumerate}
\renewcommand{\labelenumi}{(\alph{enumi})}
\item If $\phi(z) = (\epsilon_j z_{\sigma(j)})_{j \geq 1}$ for
some permutation $\sigma$ of $\NN$ and some sequence $(\epsilon_j)_{j \geq 1}$
of complex signs, then $C_\phi$ is an automorphism of
$A^+(\TT^\infty)$, and it is isometric.
\item If $C_\phi$ is an automorphism of $A^+(\TT^\infty)$ \emph{and if we moreover} assume 
that $\phi_k(z) = z_k^{d_k} u_k(z)$, with $d_k\geq 1$ and $u_k(0)\neq 0$, for each 
$k \in \NN$ and each $z \in {\textbf B}$, then $\phi(z) = (\epsilon_j z_{j})_{j \geq 1}$ for
some sequence $(\epsilon_j)_{j \geq 1}$ of complex signs.
\end{enumerate}
\end{theoreme}

\noindent{\bf Proof.} \textit{(a)} is trivial.  For \textit{(b)}, consider the compact set 
$K =\overline \DD^\infty$, endowed with the  product topology
($K$ is nothing but the spectrum of $A^+(\TT^\infty)$); clearly, 
${\textbf B}$ is dense in $K$, and since
$\phi_j = C_\phi(z_j) \in A^+(\TT^\infty)$, $\phi_j \colon {\textbf B} \to \DD$
extends continuously to $K$, and $\phi = (\phi_j)_j$ extends
continuously to a map, still denoted by $\phi$, from $K$ to $K$, and we still can write, 
for every $k \in \NN$, $\phi_k(z) = z_k^{d_k} u_k(z)$ for each $z \in K$. 
Exactly as in the Proof of Theorem \ref{theorem33}, we can show that $\phi$
is bijective, since $K$ is the spectrum of $A^+(\TT^\infty)$.  Let
now $\psi \colon K \to K$ be the inverse map of $\phi$. Since $K$ is compact, 
$\psi$ is continuous on $K$, and so on ${\textbf B}$; it is 
then easy to see, as usual, that $\psi$ is holomorphic in ${\textbf B}$ (alternatively,  
$\psi_k=(C_\phi)^{-1}(z_k)\in A^+(\TT^\infty)$, and so is analytic in 
$\DD^\infty$, and it is clear that $\psi=(\psi_k)_k$).\par
Now, it suffices to show that $\psi$ maps ${\textbf B}$ into ${\textbf B}$; indeed, 
it will follow that $\phi$ maps ${\textbf B}$ onto ${\textbf B}$, and so   
the map $\phi$ will appear as an analytic automorphism of ${\textbf B}$ (since we
already know that $\psi=\phi^{-1}$ is analytic in ${\textbf B}$), and 
Lemma \ref{lemme34} shows that $\phi_j(z)$ has the form
$\epsilon_j \frac{z_{\sigma(j)} - a_j}{1- \overline{a_j} z_{\sigma(j)}}\cdot$ Now,
$\|\phi^n_j\|_{A^+(\TT^\infty)} = \|C_\phi(z_j^n)\|_{A^+(\TT^\infty)} \leq \|C_\phi\|$, 
and as in the Proof of Theorem \ref{theorem33}, we shall conclude that $a_j = 0$ for each 
$j$. Finally, the assumption $\phi_k(z) = z_k^{d_k} u_k(z)$ for each $k$ will imply that 
$\sigma$ is the identity map.\par\smallskip
So we have to show that $\psi({\textbf B})\subseteq {\textbf B}$. If it were not the 
case, it would exist an element $w=(w_j)_j\in {\textbf B}$ such that 
$\psi(w)\notin {\textbf B}$. Hence there would exist $\delta>0$ and an infinite subset 
$J\subseteq \NN$ such that 
\begin{equation} 
|\psi_j(w)|>\delta \textrm{ for every } j\in J.
\label{delta}
\end{equation}
Let $\delta'=\delta/\|C_\psi\|$.\par
Since $w\in {\textbf B}$, we should find an integer $N\geq 1$ such that 
\begin{displaymath}
n\geq N\hskip 3mm \Rightarrow\hskip 3mm |w_n|\leq \delta'.
\end{displaymath}
Let $\kappa=\max_{1\leq n\leq N} |w_n|$. Since $\kappa<1$, there would exist $p\geq 1$ 
such that $\kappa^p<\delta'$. Consider the finite set:
\begin{displaymath}
F=\{\alpha=(m_1,\ldots,m_N,0,\ldots)\,;\ m_1+\cdots+m_N\leq p\}.
\end{displaymath}
We assert that:
\begin{equation}
F \textrm{ intersects the spectrum of } \psi_j \textrm{ for every } j\in J.
\label{meet}
\end{equation}
Indeed, writing:
\begin{displaymath}
\psi_j(z)=\sum a_j(n_1,\ldots,n_k,0,\ldots)z_1^{n_1}\ldots z_k^{n_k},
\end{displaymath}
we have:
\begin{itemize}
\item if $\alpha=(n_1,\ldots,n_l,\ldots)$ with $l>N$ and $n_l\neq 0$, then 
$|w_l|\leq \delta'$, and so:
\begin{displaymath}
|w^\alpha|\leq |w_1^{n_1}\ldots w_l^{n_l}|\leq |w_l^{n_l}|\leq |w_l|\leq \delta';
\end{displaymath}
\item if $n_1+\cdots n_N\geq p$, then:
\begin{displaymath}
|w_1^{n_1}\cdots w_N^{n_N}|\leq \kappa^{n_1+\cdots+n_N}\leq \kappa^p<\delta'.
\end{displaymath}
\end{itemize}
Hence, in both cases, $\alpha\notin F$ implies $|w^\alpha|<\delta'$. Therefore, 
if $F$ does not intersect the spectrum of $\psi_j$, we get:
\begin{displaymath}
|\psi_j(w)|\leq \sum_{\alpha\notin F} |a_j(\alpha)|\,|w^\alpha| \leq 
\delta' \|\psi_j\|_{A^+(\TT^\infty)} \leq \delta' \|C_\psi\| =\delta
\end{displaymath}
(since $\|\psi_j\|_{A^+(\TT^\infty)}=\|C_\psi(z_j)\|_{A^+(\TT^\infty)}\leq 
\|C_\psi\|\,\|z_j\|_{A^+(\TT^\infty)}=\|C_\psi\|$), which contradicts 
(\ref{delta}).\par
To end the proof, remark now that the assumption $\phi_k(z)=z_k^{d_k}u_k(z)$ 
for every $k\in \NN$ implies that: 
\begin{displaymath}
z_k=\phi_k\big[\psi(z)\big]=[\psi_k(z)]^{d_k} u_k[\psi(z)].
\end{displaymath}
But this is impossible, since $J$ is infinite and, for $k\in J$, $\psi_k(z)$ depends on 
$(z_1,\ldots,z_N)$, and hence $\phi_k\big[\psi(z)\big]=[\psi_k(z)]^{d_k} u_k[\psi(z)]$ 
also (since $d_k\geq 1$ and $u_k(0)\neq 0$).\par 
That ends the proof of Theorem \ref{theo32}.
\hfill$\square$
\bigskip 

\noindent{\bf Remark.} We shall see later, in Section 4, Theorem \ref{auto-isom}, that 
the converse of \textit{(a)} in Theorem \ref{theo32} is true.
\bigskip 

Although Theorem \ref{theo32} is  not completely satisfactory, it will be sufficient for characterizing the 
composition automorphisms of the Wiener-Dirichlet algebra $\cA^+$. In fact, we have:

\begin{theoreme}\label{theo33} Let $C_\phi \colon \cA^+ \to \cA^+$ be a composition 
operator. Then $C_\phi$ is an automorphism of $\cA^+$ if and only if 
$\phi$ is a vertical translation: $\phi(s) = s + i \tau$, where $\tau$ is a real number.
\end{theoreme}

Note that a similar result was obtained by F. Bayart \cite{B1} for the
Hilbert space $\cH^2$ of square-summable Dirichlet series
$f(s) = \sum_1^\infty a_n n^{-s}$ such that  
$\sum_1^\infty |a_n|^2 < +\infty$, but his proof does not seem to extend to our 
setting, and our strategy for proving Theorem \ref{theo33} will
be to deduce it from Theorem \ref{theo32}, with the help of the
transfer operator $\Delta$ mentioned in the Introduction. The following Lemma 
(with the notation used in the Introduction) allows the transfer from 
composition operators on $\cA^+$ to composition operators on $A^+(\TT^\infty)$.

\begin{lemme}\label{lem32} Suppose that $C_\phi \colon \cA^+ \to \cA^+$ is a
composition operator, with $\phi(s) = c_0 s + \varphi(s)$,
$c_0 \in \NN_0$, $\varphi \in \cD$.  Let 
$T = \Delta C_\phi\Delta^{-1} \colon A^+(\TT^\infty) \to A^+(\TT^\infty)$.  Then:
\begin{enumerate}
\renewcommand{\labelenumi}{(\alph{enumi})}
\item $T = C_{\tilde \phi}$, where $\tilde \phi\colon {\textbf B} \to \DD^\infty$
is an analytic map such that $\tilde \phi(z^{[s]}) = z^{[\phi(s)]}$, for any $s \in \CC_0$.
\item If moreover $c_0 \geq 1$ (which is the case if $C_\phi$ is surjective), 
$\tilde \phi$ maps ${\textbf B}$ into ${\textbf B}$.
\end{enumerate}
\end{lemme}

\noindent{\bf Proof.} \textit{(a)} Define 
$f_k(s) = p_k^{-\phi(s)} \in \cA^+$, $\phi_k = \Delta f_k$ and
\begin{eqnarray}
\tilde \phi = (\phi_1,\phi_2,\ldots) \label{equation34}.
\end{eqnarray}
We have 
\begin{displaymath}
\tilde \phi(z^{[s]}) = (\Delta f_k (z^{[s]}))_{k \geq 1} =
(f_k(s))_{k \geq 1} = z^{[\phi(s)]}
\end{displaymath}
by (\ref{equation32}), and
$\|\phi_k\|_\infty = \|f_k\|_\infty \leq 1$ by (\ref{equation33}).
Moreover, no $\phi_k$ is constant, so the open mapping theorem
implies that $|\phi_k(z)|< 1$ for $z \in {\textbf B}$, \textit{i.e.}
$\tilde \phi(z) \in \DD^\infty$.\\  
Finally, if 
$f(z) =\sum_{n=1}^\infty a_n z_1^{\alpha_1} \ldots z_r^{\alpha_r} \in A^+(\TT^\infty)$ 
(where $n =p_1^{\alpha_1} \ldots p_r^{\alpha_r}$ is the decomposition in prime factors), we 
have the following ``diagram'':
\begin{displaymath}
f~~{\buildrel \Delta^{-1} \over \longmapsto}~~ 
\sum_{n=1}^\infty a_n n^{-s} {\buildrel C_\phi \over \longmapsto}
\sum_{n=1}^\infty a_n f_1^{\alpha_1} \ldots f_r^{\alpha_r}
{\buildrel \Delta \over \longmapsto} 
\sum_{n=1}^\infty a_n \phi_1^{\alpha_1} \ldots \phi_r^{\alpha_r} =
f \circ \tilde \phi,
\end{displaymath}
\textit{i.e.} $T(f) = C_{\tilde \phi}(f)$.

\indent\textit{(b)} First observe that $C_\varphi$ also maps $\cA^+$ into
$\cA^+$ (see the remark before Corollary \ref{prop23}). Secondly, we have
$f_k(s) = p_k^{-c_0 s} p_k^{-\varphi(s)} = p_k^{-c_0 s}g_k(s)$, with
$g_k \in \cA^+$ and 
$\|g_k\|_{\cA^+} = \|C_\varphi(p_k^{-s})\|_{\cA^+}\leq C$.  It follows that, 
for $z \in {\textbf B}$~:
$\Delta f_k(z) = z_k^{c_0} \Delta g_k(z)$, and \textit{via} 
(\ref{equation33}) that:
\begin{displaymath}
|\Delta f_k(z)| \leq |z_k|^{c_0} \|\Delta g_k\|_\infty
= |z_k|^{c_0} \|g_k\|_\infty \leq |z_k|^{c_0} \|g_k\|_{\cA^+}
\leq C|z_k|^{c_0}.
\end{displaymath}
Since $c_0  \geq 1$, we see that
$\Delta f_k(z) \to 0$ as $k \to \infty$, \textit{i.e.} $\tilde \phi(z) \in {\textbf B}$. 
Finally, whenever $C_\phi$ is surjective, $\phi : \CC_0 \to \CC_0$ is injective: indeed, 
$\cA^+$ separates the points of $\CC_0$
($2^{-a} = 2^{-b}$ and $3^{-a} = 3^{-b}$ imply $a = b$, since
$\log 2/\log 3$ is irrational), and we can argue as in Theorem \ref{theorem33}.\par
To end the proof of Lemma \ref{lem32}, it remains to remark that if $c_0= 0$, $\phi$ is never 
injective on $\CC_0$, according to well-known results on the theory of analytic, almost-periodic
functions (see \textit{e.g.} Favard \cite[p. 13]{Fa}). Therefore, we have
$c_0 \geq 1$ if $C_\phi$ is surjective.\hfill $\square$
\bigskip

\noindent{\bf Proof of Theorem \ref{theo33}}. The sufficiency of the condition 
is trivial. Conversely, if $C_\phi$ is an automorphism of $\cA^+$, let 
$C_{\tilde \phi} = \Delta C_\phi \Delta^{-1}$, as in Lemma \ref{lem32}.  Since $C_\phi$ is
surjective, we know from Lemma \ref{lem32} that $\tilde \phi$ maps ${\textbf B}$ into
${\textbf B}$; we can apply Theorem \ref{theo32}, because $C_{\tilde \phi}$ is 
an automorphism of $A^+(\TT^\infty)$ onto itself and moreover 
$\tilde \phi_k(z) =\Delta f_k(z) = z^{c_0}_k \Delta g_k(z)$, with $c_0\geq 1$ 
(by Lemma \ref{lem32} again) and 
\begin{displaymath}
\Delta g_k(0)=\lim_{{\cR e}\,s\to +\infty} g_k(s)=
\lim_{{\cR e}\,s\to +\infty} p_k^{-\varphi(s)}= p_k^{-c_1}\neq 0. 
\end{displaymath}
We conclude that:
\begin{eqnarray}
\tilde \phi(z) = (\epsilon_1 z_1, \ldots, \epsilon_n z_n, \ldots),
\label{equation35}
\end{eqnarray}
for some sequence of signs $(\epsilon_n)_n$, where $z = (z_1,\ldots,z_n,\ldots)$. \par
If we now test this equality at the points
$z^{[s]} = (p_j^{-s})_j$, $s \in \CC_0$, and use (\ref{equation32}),
we see that
\begin{eqnarray}
p_j^{-\phi(s)} = \epsilon_j p_j^{-s},\hskip 3mm s \in \CC_0,\quad j \in \NN.
\label{equation36}
\end{eqnarray}
Taking the moduli in (\ref{equation36}), we get
${\cR e}\,\phi(s) = {\cR e}\,s$.  Since $\phi(s) - s$ is analytic on the
domain $\CC_0$, this implies $\phi(s) -s = i\tau$, with 
$\tau \in \RR$, thus ending  the Proof of Theorem \ref{theo33}.\hfill$\square$

\section{Isometries of $A^+(\hbox{\large $\TT\,$}^k)$, 
$A^+(\hbox{\large $\TT\,$}^\infty), \cA^+$} 

In this section, we shall characterize the composition operators
which are isometric on $A^+(\TT^k)$ and then those which are isometric on 
$A^+(\TT^\infty)$ (under an additional assumption) and on $\cA^+$. If 
$f(z) = \sum a_\alpha z^\alpha \in A^+(\TT^k)$, it
will be convenient to note $a_\alpha = \widehat f(\alpha)$. The
spectrum of $f$ (denoted by $Sp\,f$) is the set of $\alpha$'s such
that $\widehat f(\alpha) \not= 0$.  $e$ will denote the point
$(1,\ldots,1)$ of  $\overline{\DD}^k\!\!$. An elaboration of the
method of Harzallah \cite{K} allows us to show:

\begin{theoreme} \label{theo41}
Assume that $\phi = (\phi_j)_j \colon \DD^k \to \overline\DD^k$, induces a
composition operator $C_\phi \colon A^+(\TT^k) \to A^+(\TT^k)$. Then 
$C_\phi \colon A^+(\TT^k) \to A^+(\TT^k)$ is an isometry if and only if 
there exists a square matrix $A = (a_{ij})_{1 \leq i, j \leq k}$, with 
$a_{ij} \in \NN_0$ and $\det A \not=0$, and complex signs 
$\epsilon_1,\ldots,\epsilon_k$ such that:
\begin{eqnarray}
\phi_i(z) = \epsilon_i z_1^{a_{i1}} \ldots z_k^{a_{ik}},\hskip 3mm 
1 \leq i \leq k,\quad z = (z_1,\ldots,z_k) \in \DD^k. \label{equation41}
\end{eqnarray}
\end{theoreme}

To prove this theorem, it will be convenient to use the following two 
Lemmas.

\begin{lemme} \label{lemma42}
$C_\phi$ is an isometry if and only if:
\begin{enumerate}
\renewcommand{\labelenumi}{(\alph{enumi})}
\item $\phi_i = \epsilon_i F_i$, $1 \leq i \leq k$, where $\epsilon_i$ 
is a complex sign, $\widehat{F_i} \geq 0$, and 
$F_i(e) = \|F_i\|_\infty = 1$;
\item if $\alpha,\alpha^\prime \in \NN^k_0$ are distinct, the spectra of 
$\phi^\alpha$ and $\phi^{\alpha^\prime}$ are disjoint.
\end{enumerate}
\end{lemme}

\noindent{\bf Proof.}  Suppose that \textit{(a)} and \textit{(b)} hold, and let
$f(z) = \sum \hat f(\alpha) z^\alpha \in A^+(\TT^k)$.  We have by
\textit{(b)}:
\begin{displaymath}
\|C_\phi f\|_{A^+(\TT^k)} = \sum |\hat f(\alpha)|\, \|\phi^\alpha\|_{A^+(\TT^k)} 
= \sum |\hat f(\alpha)|\,\|F^\alpha\|_{A^+(\TT^k)},
\end{displaymath}
since, with the obvious notation, 
$\phi^\alpha = \epsilon^\alpha F^\alpha$. Since $\widehat{F^\alpha} \geq 0$, we have,
using \textit{(a)}~:
\begin{displaymath}
\|F^\alpha\|_{A^+(\TT^k)} = F^\alpha(e) = 1,
\end{displaymath}
so that:
\begin{displaymath}
\|C_\phi f\|_{A^+(\TT^k)} = \sum |\hat f(\alpha)| = \|f\|_{A^+(\TT^k)}.
\end{displaymath}

Conversely, suppose that $C_\phi$ is an isometry. For each $i \in [1,k]$
and each $n \in \NN$, we have $\|\phi^n_i\|_{A^+(\TT^k)} = \|z^n_i\|_{A^+(\TT^k)} = 1$,
whence
$\|\phi_i\|_\infty = 
\lim\limits_{n \to \infty} \|\phi^n_i\|_{A^+(\TT^k)}^{1/n}= 1$, by the spectral 
radius formula.  Since $\|\phi_i\|_\infty \leq \|\phi_i\|_{A^+(\TT^k)} = 1$, the only 
possibility is that 
$\phi_i = \epsilon_i F_i$, with $|\epsilon_i| = 1$, $\widehat{F_i} \geq 0$, and 
$\|\phi_i\|_{A^+(\TT^k)}= 1 = \|F_i\|_{A^+(\TT^k)} = F_i(e)$.  Therefore, 
\textit{(a)} holds.  Now suppose that we
can find $\alpha \not= \alpha^\prime$ such that
$Sp\, \phi^\alpha \cap Sp\, \phi^{\alpha^\prime}$
contains an element $\beta_0 \in \NN^k_0$, and set
$\rho = \widehat{\phi^\alpha}(\beta_0)$, 
$\rho^\prime =\widehat{\phi^{\alpha^\prime}}(\beta_0)$.  Without loss of
generality, we may assume that $|\rho| \geq |\rho^\prime|$.  Let $\theta$
be a complex sign such that 
$|\rho + \theta \rho^\prime| = |\rho| -|\rho^\prime|$.  Then, we have 
$\|z^\alpha + \theta z^{\alpha^\prime}\|_{A^+(\TT^k)} = 2$, whereas
\begin{align*}
\|C_\phi(z^\alpha + \theta z^{\alpha^\prime})\|_{A^+(\TT^k)}
& = \|\phi^\alpha + \theta \phi^{\alpha^\prime}\|_{A^+(\TT^k)} \cr
&=\sum_{\beta \not= \beta_0} |\widehat{\phi^\alpha}(\beta) + \theta
\widehat{\phi^{\alpha^\prime}}(\beta)| + |\rho + \theta
\rho^\prime| \cr
& \leq  \sum_{\beta \not= \beta_0} |\widehat{\phi^\alpha}(\beta)| +
\sum_{\beta \not= \beta_0} |\widehat{\phi^{\alpha^{\prime}}}(\beta)|
+ |\rho| - |\rho^\prime| \cr 
& =  1 - |\rho| + 1 - |\rho^\prime| + |\rho|
- |\rho^\prime| = 2(1 - |\rho^\prime|) < 2,
\end{align*}
contradicting the isometric character of $C_\phi$.\hfill $\square$

\begin{lemme} \label{lemme17}
If $\phi = (\phi_i)_i$ and if one of the $\phi_i$'s is not a monomial,
then we can find a pair of distinct elements $\alpha$, 
$\alpha^\prime \in \NN^k_0$ such that the spectra of $\phi^\alpha$
and $\phi^{\alpha^\prime}$ intersect.
\end{lemme}

\noindent{\bf Proof.} To avoid awkward notation, we will assume that $k = 3$, but it
will be clear that the reasoning works for any
value of $k$. Since only the spectra of the $\phi_i$'s are involved,
we can assume without loss of generality that we have:
\begin{eqnarray*}
\phi_1(z) &=& z_1^{s_1} z_2^{s_2}z_3^{s_3} + z_1^{t_1}
z_2^{t_2}z_3^{t^3},\hskip 3mm 
\mbox{\rm with $(s_1,s_2,s_3) \not= (t_1,t_2,t_3)$},
\\
\phi_2(z) &=& z_1^{u_1}z_2^{u_2}z_3^{u_3},\\
\phi_3(z) &=& z_1^{v_1} z_2^{v_2} z_3^{v_3}
\end{eqnarray*}
(in short, $\phi_1(z) = z^s + z^t$; $\phi_2(z) = z^u$; $\phi_3(z) = z^v)$. \par
If $\alpha = (a,b,c)$, the spectrum of $\phi^\alpha = (z^s + z^t)^a
z^{bu}z^{cv}$ consists of the triples
\begin{displaymath}
\rho s_j + (a-\rho)t_j + bu_j + cv_j = \rho(s_j - t_j) + at_j +
bu_j + cv_j,
\end{displaymath}
with $j = 1,2,3$ and $0 \leq \rho \leq a$. Therefore, if
$\alpha^\prime = (a',b',c')$, the spectra of $\phi^\alpha$ and $\phi^{\alpha^\prime}$
will intersect if and only if we can find $0 \leq \rho \leq a$ and 
$0 \leq \rho^\prime \leq a'$ such that: 
\begin{displaymath}
\rho(s_j - t_j) + at_j + bu_j + cv_j = 
\rho^\prime(s_j - t_j) + a't_j + b'u_j + c'v_j,\hskip 3mm j = 1,2,3,
\end{displaymath}
or equivalently:
\begin{eqnarray}
(\rho - \rho^\prime)(s_j - t_j) + (a-a')t_j + (b-b')u_j =
(c'-c)v_j,\hskip 3mm j = 1,2,3. \label{equation42}
\end{eqnarray}
In (\ref{equation42}), we can drop the conditions $\rho \leq a$,
$\rho^\prime \leq a^\prime$, since we can always replace $a$ and
$a^\prime$ by $a + N$ and $a^\prime + N$, where $N$ is a large integer, without affecting the
result.  Now, let $M$ be the matrix:
\begin{displaymath}
M = \left[  \begin{array}{ccc}
s_1 - t_1 & t_1 & u_1 \\
s_2 - t_2 & t_2 & u_2 \\
s_3 - t_3 & t_3 & u_3 \end{array} \right].
\end{displaymath}
To solve equation (\ref{equation42}), we distinguish two cases.

\bigskip
\noindent{\bf Case 1 :} $\det M = 0$.
\par
We decide then to take $c' = c$.  Since the field  $\QQ$ of rational
numbers is the quotient field of $\ZZ$, we can find $\lambda,\mu,\nu \in\ZZ$, not all 
zero, such that:
\begin{displaymath}
\lambda(s_j -t_j) + \mu t_j + \nu u_j = 0,\hskip 3mm j = 1,2,3.
\end{displaymath}
If $\mu$ and $\nu$ are both zero, then $\lambda = 0$, since $s_j - t_j \not= 0$
for some $j$.  Therefore, we may assume for example that $\mu \not= 0$, and write 
$\lambda = \rho - \rho^\prime$, $\mu = a - a^\prime$, $\nu = b - b^\prime$, with 
$\alpha = (a,b,c) \in\NN^3_0$, 
$\alpha^\prime = (a^\prime,b^\prime,c^\prime) \in\NN^3_0$, and 
$\alpha \not= \alpha^\prime$ since $a \not=a^\prime$.  By construction, we have 
(\ref{equation42}), so that the spectra of $\phi^\alpha$ and $\phi^{\alpha^\prime}$ are not
disjoint.
\goodbreak

\noindent{\bf Case 2 : } $\det M \not= 0$.
\par We can then find rational numbers $q,r,s$ such that:
\begin{displaymath}
q(s_j - t_j) + rt_j + su_j = v_j,\hskip 3mm j = 1,2,3,
\end{displaymath}
and we can write $q = \frac{\lambda}{N}$, $r = \frac{\mu}{N}$, 
$s =\frac{\nu}{N}$, where $\lambda,\mu,\nu \in \ZZ$ and where $N$ is a
positive integer.  Therefore, we have:
\begin{displaymath}
\lambda(s_j-t_j) + \mu t_j + \nu u_j = Nv_j, \hskip 3mm 1 \leq j \leq 3,
\end{displaymath}
and writing $\lambda = \rho - \rho^\prime$, $\mu = a - a^\prime$,
$\nu = b-b^\prime$, $c = 0$, $c^\prime = N$, we get
(\ref{equation42}) with distinct triples $\alpha = (a,b,c)$ and
$\alpha^\prime = (a^\prime,b^\prime,c^\prime)$ of non-negative
integers. Once again, the spectra of $\phi^\alpha$ and $\phi^{\alpha^\prime}$
are not disjoint.\hfill$\square$

\bigskip
\noindent{\bf Proof of Theorem \ref{theo41}}. If the condition holds, 
$C_\phi$ is an isometry by Lemma \ref{lemma42}.\par
Conversely, suppose that $C_\phi$ is an isometry.  Then, by Lemma \ref{lemma42}, 
the spectra of $\phi^\alpha$ and $\phi^{\alpha^\prime}$ are disjoint if 
$\alpha \not=\alpha^\prime$, and by Lemma \ref{lemme17} each $\phi_i$ is a
monomial, necessarily of the form (\ref{equation41}) by \textit{(a)} of
Lemma \ref{lemma42}. Finally, if we denote by $A$ the square
matrix $(a_{ij})$, by $A^\ast = (a_{ji})$ its adjoint matrix, and
if we let $A$, $A^\ast$ act on $\ZZ^k$ by the formulas:
\begin{eqnarray}
A(\alpha) = \beta\,,\hskip 5mm A^\ast(\alpha) =\gamma, \label{equation43}
\end{eqnarray}
with $\beta_i = \sum_{j=1}^k a_{ij} \alpha_j\,$ and 
$\gamma_j = \sum_{i=1}^k a_{ij} \alpha_i$, we see that:
\begin{eqnarray}
C_\phi(z^\alpha) = \phi^\alpha = 
\epsilon^\alpha z^{A^\ast(\alpha)} \label{equation44}.
\end{eqnarray}
In fact,
\begin{displaymath}
C_\phi(z^\alpha) = \prod\limits_i \phi_i^{\alpha_i} = \prod\limits_i
\epsilon_i^{\alpha_i} \Big(\prod\limits_j z_j^{a_{ij}}\Big)^{\alpha_i}
= \epsilon^\alpha \prod\limits_j z_j^{\gamma_j}.
\end{displaymath}
Now, by Lemma \ref{lemma42}, the $\phi^{\alpha}$'s have disjoint
spectra, so that the $A^\ast(\alpha)$'s are distinct, implying
$\det A \not= 0$.\hfill$\square$
\bigskip

If we now turn to the case of $A^+(\TT^\infty)$, Lemma
\ref{lemma42} clearly still holds, but Lemma \ref{lemme17} no
longer holds~: for example, if $I_1,\ldots,I_n,\ldots$ are
disjoint subsets of $\NN$, $c_{ij}$ positive numbers such that
$\sum\limits_{j\in I_i} c_{ij} = 1$, $i=1,2,\ldots$ and if the map
$\tilde \phi$ is defined by:
\begin{eqnarray}
\tilde \phi = (\phi_i)_i,\hskip 3mm
\mbox{where}\ \phi_i(z) = \sum_{j \in I_i} c_{ij} z_j,
\label{equation45}
\end{eqnarray}
then $C_{\tilde \phi}$ is an isometry by Lemma \ref{lemma42} and
yet no $\phi_i$ is a monomial if each $I_i$ has more than one element.  
We have however a weaker result:
\goodbreak

\begin{theoreme} \label{theorem44}
Let $\phi \colon \overline \DD^\infty \to \overline \DD^\infty$ be a map 
inducing a composition operator 
$C_\phi \colon A^+(\TT^\infty) \to A^+(\TT^\infty)$, and such that moreover 
$\phi(\TT^\infty) \subseteq \TT^\infty$. Then:
\begin{enumerate}
\renewcommand{\labelenumi}{(\alph{enumi})}
\item There exists a matrix 
$A = (a_{ij})_{i,j \geq 1}$, with $a_{ij} \in \NN_0$ and $\sum_j a_{ij} < \infty$ for 
each $i$, and complex signs $\epsilon_i$ such that $\phi = (\phi_i)_i$ and
\begin{eqnarray}
\phi_i(z) = \epsilon_i \prod^\infty_{j=1} z_j^{a_{ij}},\hskip 3mm i = 1,2,\ldots
\label{equation46}
\end{eqnarray}
\item $C_\phi$ is an isometry if and only if 
$A^\ast = (a_{ji})$, acting on $\ZZ^{(\infty)}$ as in (\ref{equation43}), is injective.
\end{enumerate}
\end{theoreme}

\noindent{\bf Proof.}
\textit{(a)} If we apply Lemma \ref{lemme8}  to the (connected) group 
$\Gamma = \TT^\infty$ and its dual $G =\ZZ^{(\infty)}$, we see that for each 
$i \in \NN$ there exists
$L_i = (a_{i1},a_{i2},\ldots) \in \ZZ^{(\infty)}$, necessarily in
$\NN_0^{(\infty)}$, and a complex sign $\epsilon_i$ such that, for
each $z \in \TT^\infty$, we have:
\begin{displaymath}
\phi_i(z) = \epsilon_i < L_i,z >\, = \epsilon_i \prod_j z_j^{a_{ij}} 
\end{displaymath}
(note that, for $n \in \NN$, setting $C = \|C_\phi\|$, we have
\begin{displaymath}
\|\phi^n_i\|_{A^+(\TT^\infty)} = \|C_\phi(z_i^n)\|_{A^+(\TT^\infty)} \leq C, 
\end{displaymath}
and also, since $|\phi_i(e^{it})| = 1$~: 
\begin{displaymath}
\|\phi_i^{-n}\|_{A^+(\TT^\infty)} = \|\overline \phi_i^n\|_{A^+(\TT^\infty)}
= \|\phi_i^n\|_{A^+(\TT^\infty)} \leq C).
\end{displaymath}
This proves (\ref{equation46}).\par

\textit{(b)} We know from (\ref{equation44}) (which clearly still holds
for $k = \infty$) that 
$\phi^\alpha = \epsilon^\alpha z^{A^\ast(\alpha)}$, and we know from 
Lemma \ref{lemma42} that $C_\phi$ is an isometry if and only if the spectra of the 
$\phi^\alpha$'s are disjoint.  This gives the result.\hfill$\square$
\bigskip

We shall prove here the announced converse of part \textit{(a)} of Theorem \ref{theo32}.\par

\begin{theoreme} \label{auto-isom}
Let $\phi=(\phi_j)_j\colon \textbf{B}\to \textbf{B}$ be an analytic function which induces 
a composition operator $C_\phi$ on $A^+(\TT^\infty)$. If $C_\phi$ is an isometric  
automorphism of $A^+(\TT^\infty)$, then $\phi(z)=(\epsilon_j z_{\sigma(j)})_j$, for 
some permutation $\sigma$ of $\NN$ and some some sequence $(\epsilon_j)_{j\geq1}$ of 
complex signs.
\end{theoreme}

\noindent{\bf Proof.} It suffices to look at the proof of Theorem \ref{theo32}, 
\textit{(b)}: as in that proof, and with the same notation, it suffices to show that 
$\psi(\textbf{B})\subseteq \textbf{B}$; but if it is not the case, it follows from 
(\ref{meet}), since the set $J$ is infinite, that there exist at least two distinct 
integers $j_1,j_2\in J$ such that the spectra of $\phi_{j_1}$ and 
$\phi_{j_2}$ are not disjoint. By Lemma \ref{lemma42}, this contradicts the isometric 
nature of $C_\phi$.\hfill$\square$
\bigskip

\noindent{\bf Remark.} It is easy to see that the composition operator $C_{\tilde \phi}$ on 
$A^+(\TT^\infty)$ given by (\ref{equation45}) does not correspond in general to a
$C_\phi \colon \cA^+ \to \cA^+$.\par
For example, if 
\begin{equation}
\phi_i(z) = \frac{z_{2i-1} + z_{2i}}{2}\raise0,5mm\hbox{,} \hskip 3mm 
i=1,2,\ldots,\label{exemple} 
\end{equation}
the equation $\tilde \phi(z^{[s]}) = z^{[\phi(s)]}$
would give:
\begin{displaymath}
\frac{p^{-s}_{2i-1} + p_{2i}^{-s}}{2} = p_i^{- \phi(s)}, i = 1,2,\ldots;
\end{displaymath}
taking equivalents of both members as $s \to \infty$ would give that
\begin{displaymath}
\frac{\phi(s)}{s} \mathop{\longrightarrow}\limits_{s\to+\infty}
\frac{\log p_{2i-1}}{\log p_i}\raise0,4mm\hbox{,}
\end{displaymath} 
and it is impossible to have that, even for one $i$, since 
$\frac{\phi(s)}{s} \to c_0 \in \NN_0$~! \par
On the other hand, the additional assumption made in Theorem \ref{theorem44} does not 
allow to use the Bohr's transfer operator $\Delta$ to characterize the isometric 
composition operators on $\cA^+$. Nevertheless, we have:\goodbreak

\begin{theoreme}\label{Dirichlet-isometries}
Let $\phi\colon \CC_0\to \CC_0$ inducing a composition operator 
$C_\phi\colon \cA^+\to \cA^+$. Then $C_\phi$ is an isometry if and only if 
$\phi(s)= c_0 s +i\tau$, with $c_0\in \NN$ and $\tau\in \RR$.
\end{theoreme}

\noindent{\bf Proof.} One direction is trivial. For the other, let us 
introduce the following notation: if $f(s)=\sum_{k=1}^\infty a_k k^{-s}\in \cA^+$, 
denote by $Sp\,f$ (the spectrum of $f$) the set of indices $k$ such that $a_k\neq 0$. 
Now, the technique of the proof of Lemma \ref{lemma42} clearly works to show that:
\begin{equation}
\mbox{\textit{If $m$ and $n$ are distinct integers, the spectra of $m^{-\phi}$ 
and $n^{-\phi}$ are disjoint}}\label{disjoints}
\end{equation}
This automatically implies $c_0\neq 0$, since, otherwise, the integer $1$ would belong  
to the spectra of all the $n^{-\phi}$'s. Suppose now that $\phi$ is not of the form 
$c_0 s+ c_1$, and write:
\begin{displaymath}
\phi(s)= c_0 s + c_1 + \omega(s),
\end{displaymath}
with: 
\begin{displaymath}
\omega(s)=c_r r^{-s} +c_{r+1} (r+1)^{-s}+\cdots,\hskip 3mm r\geq 2,\quad c_r\neq 0.
\end{displaymath}
Then:
\begin{eqnarray*}
n^{-\phi(s)}&=& (n^{c_0})^{-s} n^{-c_1} \exp\big(-\omega(s)\log n\big)\\
&=&(n^{c_0})^{-s} n^{-c_1} \Big[1+\sum_{k=1}^\infty 
\frac{(-\log n)^k}{k!}\big(\omega(s)\big)^k\Big]\\
&=& (n^{c_0})^{-s} n^{-c_1} \Big[1+\cdots+\sum_{k=1}^\infty \frac{(-\log n)^k}{k!}
(c_r r^{-s}+\cdots)^k\Big].
\end{eqnarray*}
For ${\cR e}\,s$ large enough, all the series involved will be absolutely convergent; therefore 
the Dirichlet series of $n^{-\phi}$ will be obtained by expanding 
$(c_r r^{-s}+\cdots)^k$ and grouping terms. In particular, the coefficient $\lambda_n$ of 
$n^{c_0}r^{c_0}$ in $n^{-\phi}$ can be obtained only by expanding $(c_r r^{-s}+\cdots)^k$ for 
$k=1,\ldots,c_0$, so that $\lambda_n=P(\log n)$, where $P$ is a non-zero polynomial. This 
implies that, for large $n$, $\lambda_n\neq 0$, and $(nr)^{c_0}\in Sp\,n^{-\phi}$. Moreover, 
it is clear that $l^{c_0}\in Sp\,l^{-\phi}$ for every positive integer $l$. Hence 
$(nr)^{c_0}\in Sp\,n^{-\phi}\cap Sp\,(nr)^{-\phi}$ for large $n$, which contradicts 
(\ref{disjoints}).\par
Therefore $\phi(s)=c_0 s+c_1$, and $c_1$ clearly has to be purely imaginary if $C_\phi$ is 
an isometry.\hfill$\square$
\bigskip

\section{Concluding remarks and questions}

Proposition \ref{prop23} does not answer, in general, the
natural question~: if $C_\phi$ maps $\cA^+$ into $\cA^+$, is it true
that $\phi(s) = c_0s +\sum_1^\infty c_n n^{-s}$, with
$\sum_1^\infty |c_n| < \infty$~? \par
Proposition \ref{prop24} does not apply to the case of 
complex coefficients $c_r$, $c_{r^2}$. Here, recent estimates due to
Rusev \cite{Ru} might help. \par
The estimate $\|\phi^n\|_{A^+} \geq \delta \sqrt n$ of Lemma \ref{lemme31} is best possible. 
In fact (see \cite[p. 76]{K}) it is fairly easy to see that
$\|\phi^n\|_{A^+} \leq C \sqrt n$ if $\phi = e^{ig}$ and $g$ is $\cC^\infty$
(say), and a similar computation in dimension $k$ (\textit{i.e.} if we work with
$A^+(\TT^k)$) easily gives the estimate $\|\phi^n\|_{A^+(\TT^k)} \leq C_k n^{k/2}$
if $\phi = e^{ig}$ and $g$ is $\cC^\infty$.  It would be interesting to know
whether the converse holds, \textit{i.e.} if we have the following
quantitative version of Lemma \ref{lemme31}~: if $\phi
= e^{ig}$, where $g$ is a $\cC^\infty$, non-affine, real function,
then $\|\phi^n\|_{A^+} \geq \delta n^{1/2}$~? \par
In the proof of Theorem \ref{theo33}, we used the fact that an 
analytic, almost-periodic, function defined on a vertical 
half-plane is never injective, to show that $c_0 > 0$, and 
therefore that the assumption \textit{(b)} in Theorem \ref{theo32} naturally
holds.  This raises two questions:

\begin{enumerate}
\renewcommand{\labelenumi}{\alph{enumi})}
\item Can an almost-periodic function defined only on a vertical
line be injective,  \textit{i.e.} can an almost-periodic function $f \colon \RR \to \CC$
be injective? (of course, if $f$ is real-valued, this is impossible: if
$f$ is injective, it is monotonic and therefore non almost-periodic).
\item Can one, at the price of using a different Banach-Stone type
Theorem, dispense with the condition $\phi_k(z) = z_k^{d_k} u_k(z)$, with 
$d_k\geq 1$ and $u_k(0)\neq 0$ of \textit{(b)} in Theorem \ref{theo32}, 
\textit{i.e.} is the converse of \textit{(a)} in this Theorem always true?
\end{enumerate}

In view of the examples (\ref{equation45}) in Section~4, a complete description of the 
isometric composition operators $C_\phi \colon A^+(\TT^\infty) \to A^+(\TT^\infty)$ seems
out of reach. \par
We gave a proof of Theorem \ref{Dirichlet-isometries} which does not use 
Theorem \ref{theo41}. Using this theorem, we can give a variant of 
Theorem \ref{Dirichlet-isometries}: fix an integer $k \geq 1$, and denote by
$\cA^+_k$ the subalgebra of $\cA^+$ consisting of the functions
$f(s) = \sum_{P^+(n) \leq k} a_n n^{-s}$, where $P^+(n)$
denotes the largest prime factor of $n$. Equivalently, $f \in \cA^+_k$ if the Dirichlet 
expansion of $f$ only involves the primes $p_1,\ldots,p_k$. Define similarly the subspace 
$\cD_k$ of $\cD$. With those definitions, we can state:

\begin{theoreme}
\label{theo51} Let $\phi(s) = c_0 s + \varphi(s)$, $\varphi \in \cD_k$, induce a 
composition operator $C_\phi \colon \cA^+_k \to \cA^+_k$.  Then 
$C_\phi \colon \cA^+_k \to \cA^+_k$ is an isometry if and only if 
$\phi(s) = c_0 s + i\tau$, with $c_0 \in \NN$ and $\tau \in \RR$.
\end{theoreme}

\noindent{\bf Proof.} The sufficiency is trivial.  For the necessity, define
an isometry $\Delta \colon \cA^+_k \to \cA^+(\TT^k)$ by:
\begin{displaymath}
\Delta \Big( \sum_{n=1}^\infty a_n n^{-s} \Big) =
\sum_{n=1}^\infty a_n z_1^{\alpha_1} \ldots z_k^{\alpha_k},
\end{displaymath}
where $n = p_1^{\alpha_1} \ldots p_k^{\alpha_k}$ is the decomposition of $n$ 
in prime factors. 
Set $z^{[s]} = (p_1^{-s}, \ldots, p_k^{-s} ) \in \DD^k$ and
check that $\Delta C_\phi \Delta^{-1} = T$ is a composition
operator $C_{\tilde \phi} \colon A^+(\TT^k) \to A^+(\TT^k)$, isometric
if $C_\phi$ is isometric, and such that:
\begin{eqnarray}
\tilde \phi (z^{[s]}) = z^{[\phi(s)]} \label{equation51}.
\end{eqnarray}
We now use Theorem \ref{theo41} to conclude that 
$\tilde \phi = (\phi_1,\ldots,\phi_k)$, with 
$\phi_1(z) = \epsilon_1 z_1^{a_{11}}\ldots z_k^{a_{1k}}$, and where 
$a_{11}, \ldots, a_{1k}$ are non-negative integers. Exactly as in the Proof of 
Theorem \ref{theorem33}, we then conclude that $\phi(s) = c_0 s + i \tau$.\hfill$\square$
\bigskip

In the next Theorem, we shall see that there are few composition 
operators whose symbols preserve the boundary $i\RR$.

\begin{theoreme}
\label{theorem45}
Let $\phi \colon \CC_0 \to \CC_0$ inducing a composition operator
$C_\phi \colon \cA^+ \to \cA^+$, and such that moreover $\phi$ has a
continuous extension to $\overline \CC_0$, preserving the boundary
of $\CC_0$, \textit{i.e.} $\phi(i \RR) \subseteq i\RR$.  Then 
$\phi(s) = c_0 s + i \tau$, where $c_0 \in \NN_0$ and $\tau \in \RR$.
\end{theoreme}

\noindent{\bf Proof.} Let $\tilde \phi$ be associated with $\phi$ as in 
Theorem \ref{theorem33}. By continuity, the equation 
$\tilde \phi\big(z^{[s]}\big) = z^{[\phi(s)]}$, $s \in \CC_0$, still
holds for $s = it$, $t \in \RR$, to give 
$\tilde \phi\big(( p_j^{-it})_j\big)= 
\big(p_j^{-\phi(it)}\big)_j$, and so $\tilde \phi(\TT^\infty) \subseteq \TT^\infty$
since, by the Kronecker Approximation Theorem and the definition
of the product topology on $\TT^\infty$, the points
$(p_j^{-it})_j$, $t \in \RR$, are dense in $\TT^\infty$. Now, by
Theorem \ref{theorem44}, we have in particular
$\tilde \phi = (\phi_i)_i$, with
\begin{displaymath}
\phi_1(z) = \epsilon_1 z_1^{a_{11}} \ldots z_k^{a_{1k}},
\end{displaymath}
for some complex sign $\epsilon_1$ and some integer $k$. 
In particular, the equation $\tilde \phi (z^{[s]}) = z^{[\phi(s)]}$
implies that:
\begin{displaymath}
\epsilon_1 (p_1^{-s})^{a_{11}} \ldots (p_k^{-s})^{a_{1k}}
= p_1^{-\phi(s)}, \hskip 3mm s \in \CC_0.
\end{displaymath}
Passing to the moduli gives ${\cR e}\,\phi(s) = c\, {\cR e}\, s$,
with $c = \sum_{j=1}^k a_{1j} \frac{\log p_j}{\log p_1}\cdot$\par 
Theorefore, $\phi(s) - cs = i \tau$, $\tau \in \RR$, and we know that $c = c_0$ is 
necessarily an integer.\hfill$\square$
\bigskip

\noindent{\bf Acknowledgments.}  The authors thank J.P. Vigu\'e
and W. Kaup for fruitful discussion and information. We also thank E. Strouse for 
correcting a great number of mistakes in English (before we added others!), and the referee for a careful reading.

\bigskip

\noindent \textit{Fr\'ed\'eric Bayart}, LaBAG, Universit\'e Bordeaux 1, 351 cours de la 
Lib\'eration, 33405 Talence cedex, France --  Frederic.Bayart@math.u-bordeaux1.fr \par\smallskip 
\noindent \textit{Catherine Finet}, Institut de Math\'ematique, Universit\'e de Mons-Hainaut, 
``Le Pentagone", Avenue du Champ de Mars, 6,  7000 Mons, Belgique -- 
catherine.finet@umh.ac.be \par\smallskip
\noindent \textit{Daniel Li}, Laboratoire de Math\'ematiques de Lens, Universit\'e d'Artois, 
rue Jean Souvraz, SP18, 62307 Lens Cedex, France --  
daniel.li@euler.univ-artois.fr \par\smallskip 
\noindent \textit{Herv\'e Queff\'elec}, UFR de Math\'ematiques, Universit\'e de Lille 1, 
59655 Villeneuve d'Ascq Cedex, France --  queff@math.univ-lille1.fr \par

\end{document}